\documentclass[12pt]{article}
\usepackage{authblk}
\usepackage[utf8]{inputenc}
\usepackage{amsmath,amssymb,amsthm,xcolor,mathrsfs,enumitem,mathtools,centernot, oubraces}
\usepackage[T1]{fontenc}
\usepackage[margin=1in]{geometry}
\usepackage{bussproofs}
\usepackage{chngcntr}
\usepackage{setspace}
\usepackage{cite}

\usepackage[hidelinks, backref=page]{hyperref}
\hypersetup{
    colorlinks,
    citecolor=magenta,
    filecolor=blue,
    linkcolor=blue,
    urlcolor=blue
    }

\theoremstyle{plain}
\newtheorem{defn}{Definition}[section]

\newtheorem{thm}[defn]{Theorem}

\theoremstyle{definition}
\newtheorem{case}{Case}[defn]
\newtheorem{subcase}{Subcase}[case]
\newtheorem{subsubcase}{Subsubcase}[subcase]
\newtheorem{subsubsubcase}{Subsubsubcase}[subsubcase]

\newtheorem{rem}[defn]{Remark}

\newtheoremstyle{case}{}{}{}{}{}{:}{}{}

\newcommand{\MultSet}{\mathsf{MultiSet}}
\newcommand{\SeqComp}{\mathsf{SeqComp}}

\newcommand{\thmref}[1]{Theorem~\ref{#1}}
\newcommand{\secref}[1]{\S\ref{#1}}

\newcommand{\defref}[1]{Definition~\ref{#1}}
\newcommand{\remref}[1]{Remark~\ref{#1}}
\newcommand{\caseref}[1]{Case~\ref{#1}}
\newcommand{\subcaseref}[1]{Subcase~\ref{#1}}
\newcommand{\subsubcaseref}[1]{Subsubcase~\ref{#1}}

\makeatletter
\newcommand{\vitus@ax}[1]{\vitus@@ax$#1$}
\def\vitus@@ax$#1=>#2${\Axiom$#1\fCenter#2$}
\newcommand\vitus@inf[2]{\vitus@@inf{#1}$#2$}
\def\vitus@@inf#1$#2=>#3${%
  \ifcase#1\or
  \expandafter\UnaryInf\or
  \expandafter\BinaryInf\or
  \expandafter\TernaryInf\or
  \expandafter\QuaternaryInf\or
  \expandafter\QuinaryInf\fi
  $#2\fCenter#3$%
}

\makeatother

\DeclareRobustCommand\deduc{\mathrel{|}\joinrel\mkern-.5mu\mathrel{-}}

  \title{Rule-Elimination Theorems}
  \author{Sayantan Roy}
  \date{}
  \affil{Department\ of Mathematics\\
   Indraprastha Institute of Information Technology-Delhi,\\
  New Delhi, India.}

\begin{document}

\maketitle
\begin{abstract}
\textsc{Cut}-elimination theorems constitute one of the most important classes of theorems of proof theory. Since Gentzen's proof of the \textsc{cut}-elimination theorem for the system \textbf{LK}, several other proofs have been proposed. Even though the techniques of these proofs can be modified to sequent systems other than $\mathbf{LK}$, they are essentially of a very particular nature; each of them describes an algorithm to transform a given proof to a \textsc{cut}-free proof. However, due to its reliance on heavy syntactic arguments and case distinctions, such an algorithm makes the fundamental structure of the argument rather opaque. We, therefore, consider rules abstractly, within the framework of logical structures familiar from \textit{universal logic} in the sense of \cite{Beziau1994} and aim to clarify the essence of the so-called ``elimination theorems''. To do this, we first give a non-algorithmic proof of the \textsc{cut}-elimination theorem for the propositional fragment of $\mathbf{LK}$. From this proof, we abstract the essential features of the argument and define something called \textit{normal sequent structures relative to a particular rule}. We then prove two \textsc{rule}-elimination
theorems for these and show that one of these has a converse. Abstracting even more, we define \textit{abstract sequent structures} and show that for these structures, the corresponding version of the \textsc{rule}-elimination theorem has a converse as well. 
\end{abstract}

\maketitle
\tableofcontents

\section{Introduction}

\textsc{Cut}-elimination theorems are one of the most important classes of theorems of proof theory. Since Gentzen's proof of the \textsc{cut}-elimination theorem for the system \textbf{LK}, introduced in \cite{Gentzen1935}, several other proofs of the theorem have been proposed (see \cite{Indrzejczak2021} for several of these proofs, including the original one by Gentzen). These proofs are essentially of a very particular nature. Even though the techniques of these proofs can be modified to sequent systems other than $\mathbf{LK}$, each describes an algorithm to transform a given proof into a \textsc{cut}-free proof. While these algorithms are useful and often can be modified suitably to other sequent systems, they are not general enough. The following paragraph from \cite{Ciabattoni2006} nicely summarises the issue. 
\begin{quote}
    Checking whether a sequent calculus admits cut-elimination is often a rather tedious task. Indeed, the proof is based on case distinctions and has to be checked for all the possible combinations of the rules. This is usually done using heavy syntactic arguments based on case distinctions, written without filling in the details (note that even Gentzen did not formalize all the cases). This renders the cut-elimination process rather opaque. It is then natural to search for general criteria that a sequent calculus has to satisfy in order to admit cut-elimination. Moreover, such criteria, if given on a suitable level of abstraction, would also provide a deeper understanding of the nature of cut-elimination. 
\end{quote}

Furthermore, despite its important consequences, what is essentially done in the \textsc{cut}-elimination theorems is nothing but the \textit{elimination} of a certain rule. One may ask the same question for any other rule. For example, one may ask whether \textsc{weakening} is eliminable from $\mathbf{LK}$ (the case of \textsc{contraction} has been dealt with in \cite{Kashima1997}). One may even ask whether there exists some necessary and sufficient condition(s) for eliminating some particular rule from a particular sequent system. For questions like this, the usual algorithmic proofs of `elimination theorems' hardly shed any light. So, as we have quoted above, it is natural to search for a suitable abstract set-up in which questions like these can be posed and (hopefully) answered. And, while one is at it, why focus on the elimination (or non-elimination) of \textit{particular} rules in \textit{particular} sequent systems? This is the starting point of the genesis of this paper. More specifically, the questions in which we are interested are the following:
\begin{enumerate}[label=$\bullet$]
    \item What makes the elimination of \textsc{cut} possible in \textbf{LK}? 
    \item Is it possible to characterise sequent systems for which \textsc{cut}-elimination holds?
    \item Is it possible to give necessary and sufficient conditions of eliminating any rule from a given sequent system? What does a ‘rule’ mean? What are we supposed to understand by a ‘sequent system’?
\end{enumerate}

The motivation for these questions comes from \textit{universal logic} in the sense of \cite{Beziau1994}. One of the most important goals of universal logic, according to \cite{Beziau2012} is ``to clarify the fundamental concepts of logic and to construct general proofs.'' In this paper, we aim to clarify the essence of the so-called ``elimination theorems'' and construct general proofs of the same. We will do the former by providing proofs of the \textsc{cut}-elimination theorem for the propositional fragment of $\mathbf{LK}$ 
and the latter by proving what we call \textsc{rule}-elimination theorems. 

As we have mentioned earlier, there already are works of a similar spirit. For example, in \cite{AvronLev2005}, the authors use semantic methods to provide a necessary and sufficient condition for \textsc{cut}-elimination, although for multiple-conclusion calculi. In \cite{Ciabattoni2006}, the authors provide necessary and sufficient conditions (both syntactic and semantic) for \textsc{cut}-elimination for a certain class of ``sequent calculus''. Without going into a detailed comparison between these landmark works on the subject, we would like to point out that the work presented here differs from these two (and several others of the same spirit) in the following aspects. First, our main objective is to prove `\textsc{rule}-elimination theorems'. Second, in introducing normal sequent structures and abstract sequent structures, we have not assumed any particular properties of them like the subformula property (as, e.g., was required in \cite{Ciabattoni2006}). Third, one of the main motivations for introducing the notion of abstract sequent structures was to find the minimal set of properties that is required to admit \textsc{cut}-elimination. The notion of proof in this context is more general than the usual one. 

The rest of the paper is organised as follows. In \hyperlink{sec:cut_elim(PLK)}{\secref{sec:cut_elim(PLK)}}, we will give two new non-algorithmic proofs of the \textsc{cut}-elimination theorem for the propositional fragment of $\mathbf{LK}$. En route, we also give a new proof of the consistency of the same (see \hyperlink{thm:PLK_cons}{\thmref{thm:PLK_cons}}). From these two proofs, we abstract the essential features of the argument and define \textit{normal sequent structures}. Two \textsc{rule}-elimination theorems for these particular classes of sequent structures, and a converse of one of them, will then be proved in \S3. In the next section, we define \textit{abstract sequent structures} and point out the essential features that made the proof of the \textsc{rule}-elimination theorems for normal sequent structures work. This will pave the way towards formulating the most general version of the \textsc{rule}-elimination theorems for certain types of abstract sequent structures. We will show that the \textsc{rule}-elimination theorem also has a converse for these structures. Finally, in the concluding section, we will point out some directions for further research.

\hypertarget{sec:cut_elim(PLK)}{\section{\textsc{Cut}-elimination Theorem for the Propositional Fragment of \texorpdfstring{$\mathbf{LK}$}{LK}}{\label{sec:cut_elim(PLK)}}}

In this section, we provide two non-algorithmic proofs of the \textsc{cut}-elimination theorem for the propositional fragment of $\mathbf{LK}$, which we will  denote by  $\mathbf{PLK}$. The motivation for the first proof comes from Buss's proof of \textsc{cut}-elimination theorem for a variant of the system $\bf{PLK}$ (see pp. 13-15 of \cite{Buss1998}). That proof is also non-algorithmic. However, there are at least three differences between that proof and ours. These are listed below.
\begin{enumerate}[label=(\arabic*)]
    \item Buss's proof is `semantic'; ours is `syntactic'.
    \item Buss's proof (especially the \textsc{completeness theorem}) relies on the \textsc{inversion theorem}, neither of which we require. 
    \item Although Buss's proof works for the system he mentions in \cite{Buss1998}, it doesn't seem to work for the original system of Gentzen (see the first paragraph on page 14 of \cite{Buss1998}). Our proof works for both. 
\end{enumerate}

\subsection{Basics of \textbf{PLK}}

The presentation of the basics of \textbf{PLK} in this section closely follows \cite{Indrzejczak2021}.

In what follows, we will use a rather standard form of a propositional language for classical logic. Let $\mathscr{L}_{\bf{CPL}}$ denote an algebra of formulae:
$(\mathscr{L},\neg, \land,\lor,\to)$ with a countably infinite set of propositional symbols, which we denote as $\mathscr{P}$. 
Operations of this algebra correspond to the well-known connectives of \textit{negation}, \textit{conjunction}, \textit{disjunction} and \textit{implication} respectively. We will sometimes write $\alpha\in \mathscr{L}_{\mathbf{CPL}}$ to mean that $\alpha\in \mathscr{L}$.

The \textit{complexity} of a formula is a map $\mathsf{Comp}:\mathscr{L}\to\mathbb{N}$ as follows:
\begin{align*}\mathsf{Comp}(p)&:=0&\text{if}~p~\text{is a propositional variable}\\\mathsf{Comp}(\neg A)&:=\mathsf{Comp}(A)+1\\\mathsf{Comp}(A\circ B)&:=\mathsf{Comp}(A)+\mathsf{Comp}(B)\end{align*}where $\circ\in \{\land,\lor,\to\}$. 

The \textbf{PLK}\textit{-sequents} (or, when the context is clear, simply \textit{sequents}) are finite ordered lists (including the empty list) of the form $\Gamma\Rightarrow \Delta$ where the elements of both $\Gamma$ and $\Delta$ are elements of $\mathscr{L}$. We denote the empty list as $\Lambda$. 

The \textit{sequent complexity} of a sequent is a map $\mathsf{SeqComp}:\mathscr{L}\to\mathbb{N}$ as follows:
\begin{align*}\mathsf{SeqComp}(\Gamma\Rightarrow \Delta)&:=\displaystyle\sum_{\alpha\in \Gamma}\mathsf{Comp}(\alpha)+\displaystyle\sum_{\alpha\in \Delta}\mathsf{Comp}(\alpha)\end{align*}Notice that here $\Gamma$ and $\Delta$ both are finite ordered lists and not sets. 

The sequent system \textbf{PLK} consists of the following rules (here $p$ denotes a propositional variable and all other letters denote arbitrary formulas):

\begin{center}
        \underline{\textbf{AXIOMS}}
        $$p\Rightarrow  p$$
    \end{center}
    \begin{center}
        \underline{\textbf{STRUCTURAL RULES}}
        \begin{align*}&{\color{red}{\mathsf{WL}}}~~~~\dfrac{{\color{white}{A,}} \Gamma\Rightarrow \Delta}{A, \Gamma\Rightarrow \Delta}&\dfrac{\Gamma\Rightarrow \Delta {\color{white}{, A}}}{ \Gamma\Rightarrow \Delta, A}~~~~{\color{red}{\mathsf{WR}}}\\&{\color{red}{\mathsf{CL}}}~~~~\dfrac{A,A,\Gamma\Rightarrow \Delta}{{\color{white}{A,}} A, \Gamma\Rightarrow \Delta}&\dfrac{\Gamma\Rightarrow \Delta,A,A}{ \Gamma\Rightarrow \Delta, A {\color{white}{, A}}}~~~~{\color{red}{\mathsf{CR}}}\\&{\color{red}{\mathsf{EL}}}~~~~\dfrac{\Gamma, A, B,\Pi\Rightarrow \Delta}{\Gamma, B, A,\Pi\Rightarrow \Delta}&\dfrac{\Gamma\Rightarrow \Delta, A, B, \Lambda}{ \Gamma\Rightarrow \Delta, B,A, \Lambda}~~~~{\color{red}{\mathsf{ER}}}\end{align*}
        \underline{\textbf{LOGICAL RULES}}
        \begin{align*}&{\color{red}{\neg\mathsf{L}}}~~~~\dfrac{{\color{white}{\neg A,}} \Gamma\Rightarrow \Delta, A}{\neg A, \Gamma\Rightarrow \Delta {\color{white}{,A}}}&\dfrac{A,\Gamma\Rightarrow \Delta{\color{white}{,\neg A}}}{{\color{white}{ A,}}\Gamma\Rightarrow \Delta,\neg A }~~~~{\color{red}{\neg\mathsf{R}}}\\&{\color{red}{\land\mathsf{L}}}~~~~\dfrac{~~~~~A,\Gamma\Rightarrow \Delta}{A\land B, \Gamma\Rightarrow \Delta}&\dfrac{~~~~~B,\Gamma\Rightarrow \Delta}{A\land B, \Gamma\Rightarrow \Delta}~~~~{\color{red}{\land\mathsf{L}}}\\&{\color{red}{\land\mathsf{R}}}~~~~\dfrac{\Gamma\Rightarrow \Delta,A\quad\quad\Gamma\Rightarrow \Delta,B}{\Gamma\Rightarrow \Delta, A\land B}&\dfrac{A,\Gamma\Rightarrow \Delta\quad\quad B,\Gamma\Rightarrow \Delta}{A\lor B,\Gamma\Rightarrow \Delta}~~~~{\color{red}{\lor\mathsf{L}}}\\&{\color{red}{\lor\mathsf{R}}}~~~~\dfrac{\Gamma\Rightarrow \Delta, A~{\color{white}{\lor B}}}{\Gamma\Rightarrow \Delta, A\lor B}&\dfrac{\Gamma\Rightarrow \Delta, B~{\color{white}{\lor A}}}{\Gamma\Rightarrow \Delta, A\lor B}~~~~{\color{red}{\lor\mathsf{R}}}\\&{\color{red}{\to\mathsf{L}}}~~~~\dfrac{\Gamma\Rightarrow \Delta,A\quad\quad B,\Gamma\Rightarrow \Delta}{A\to B,\Gamma\Rightarrow \Delta}&\dfrac{A,\Gamma\Rightarrow \Delta, B{\color{white}{,A\to B}}}{{\color{white}{A,}}\Gamma\Rightarrow \Delta, A\to B{\color{white}{,B}}}~~~~{\color{red}{\to\mathsf{R}}}\end{align*}
        \underline{\textbf{CUT RULE}}
        \begin{align*}\dfrac{\Gamma\Rightarrow \Delta, A\quad\quad A,\Gamma\Rightarrow \Delta}{\Gamma\Rightarrow \Delta}\end{align*}
    \end{center}
We use the familiar convention of denoting the upper sequents as the \textit{premise(s)} and the lower sequent as the \textit{conclusion} of the rule under consideration. 

A formula introduced by the application of a logical rule is the \textit{principal formula} of this rule application. Formulae used for the derivation of the principal formula are \textit{side formulae} of this rule application. The principal and side formulae are the \textit{active formulae} of this rule application. By convention, we take the formula which occurs on both sides of any axiomatic sequent as an active formula.

\textit{Proofs} in \textbf{PLK} are built as trees of sequents (i.e. with nodes labelled with sequents) with axioms as leaves and the proven sequent as a root. Formally, this notion can be defined as follows.
\begin{enumerate} 
\item An axiomatic sequent $S$ is a proof of a sequent $S$.
\item If $\mathcal{D}$ is a proof of a sequent $S$, then $\dfrac{\mathcal{D}}{S'}$ is a proof of a sequent $S'$
, provided that $S$ is an instance of the premise and $S'$ an instance of the conclusion of some one-premise rule.
\item If $\mathcal{D}$ is a proof of a sequent $S$ and $\mathcal{D'}$ is a proof of a sequent $S'$, then
$\dfrac{\mathcal{D}~~~\mathcal{D'}}{S''}$ is a proof of a sequent $S''$, provided that $S$ and $S'$ are instances of
the premises and $S''$ an instance of the conclusion of some two-premise rule.
\end{enumerate}

Admitting any sequents (not only axiomatic ones) as leaves of a tree, we say that there is a \textit{derivation} of $S$ from $S_1,\ldots S_n$ if there is a tree where some leaves are decorated with sequents that are not axioms but belong to the list $S_1,\ldots S_n$. If $n = 0$, we have a proof of $S$ as a special case of a derivation. 

\begin{rem}One may ask the reason for taking the `additive' version of the \textsc{cut}-rule (and of ${\color{red}{\to\mathsf{L}}}$) (see the discussion on p. 65 of \cite{Indrzejczak2021}) instead of the `multiplicative' one (which was the original version considered by Gentzen), i.e., 
\begin{prooftree}

\def\fCenter{\Rightarrow}
\AxiomC{$\Gamma\fCenter \Delta, A$}
\AxiomC{$A, \Gamma'\fCenter \Delta'$}
\BinaryInf$\Gamma,\Gamma'\fCenter\Delta,\Delta'$
\end{prooftree}
The reason is mainly due to certain technical conveniences (a detailed account is given in \hyperlink{rem:additive_choice}{\remref{rem:additive_choice}}) and also because the multiplicative version of \textsc{cut} in \textbf{LK} is eliminable iff the additive version of \textsc{cut} in \textbf{LK} is eliminable. For the sake of completeness, we give an outline of the proof of this claim. 

For this, let us denote the additive version of \textsc{cut}, $A\textsc{cut}$ and the multiplicative one as $M\textsc{cut}$. Let us also denote the system of $\mathbf{PLK}$ with $A\textsc{cut}$ as its ``\textsc{cut}'' rule as $\mathbf{PLK}^A$ and the one with $M\textsc{cut}$ as $\mathbf{PLK}^M$.

So, suppose $A\textsc{cut}$ is eliminable in $\mathbf{PLK}^A$, but $M\textsc{cut}$ is not. This implies that there exists a proof which can't be done without using $M\textsc{cut}$. This implies, in particular, that there exists a proof whose last two lines look like the following:
\begin{prooftree}
\def\fCenter{\Rightarrow}
\AxiomC{$\Gamma\fCenter \Delta, A$}
\AxiomC{$A, \Gamma'\fCenter \Delta'$}
\BinaryInf$\Gamma,\Gamma'\fCenter\Delta,\Delta'$
\end{prooftree}We transform the proof in the following way. 
\begin{prooftree}
\def\fCenter{\Rightarrow}
\AxiomC{$\Gamma\fCenter \Delta, A$}
\LeftLabel{{\color{red}{$\textsc{\normalsize{structural}}\atop\textsc{\normalsize{rules}}$}\quad}}
\UnaryInfC{$\Gamma,\Gamma'\fCenter \Delta,\Delta', A$}
\AxiomC{$A, \Gamma'\fCenter \Delta'$}
\RightLabel{{\color{red}{\quad$\textsc{\normalsize{structural}}\atop\textsc{\normalsize{rules}}$}}}
\UnaryInfC{$A,\Gamma,\Gamma'\fCenter \Delta,\Delta'$}
\BinaryInf$\Gamma,\Gamma'\fCenter\Delta,\Delta'$
\end{prooftree}Since $A\textsc{cut}$ is eliminable in $\mathbf{PLK}^A$, it follows that there exists an $A\textsc{cut}$-free $\mathbf{PLK}^A$-proof of $\Gamma,\Gamma'\Rightarrow\Delta,\Delta'$. This means that there exists a $\mathbf{PLK}^A$-proof of $\Gamma,\Gamma'\Rightarrow\Delta,\Delta'$ using some or all of \textsc{weakening}, \textsc{exchange}, \textsc{contraction} and the logical rules. Consequently, it follows that there exists a $M\textsc{cut}$-free $\mathbf{PLK}^A$-proof of $\Gamma,\Gamma'\Rightarrow\Delta,\Delta'$ as well, a contradiction. This contradiction establishes our claim. 
\end{rem}

Throughout this section, we will assume that we are working in $\mathbf{PLK}^A$, which we will continue to denote as $\mathbf{PLK}$.

\subsection{First Proof of the \textsc{Cut}-elimination Theorem for \textbf{PLK}}

Henceforth, $\emptyset\vdash_{\textbf{PLK}}\Gamma\Rightarrow \Delta$ will be used to mean that $\Gamma\Rightarrow \Delta$ has a \textbf{PLK}-proof. Our first proof of the \textsc{cut}-elimination theorem for $\mathbf{PLK}$ needs the use of the following theorems. The proofs of these results, being straightforward, are omitted. 

\hypertarget{ml1}{\begin{thm}{\label{ml1}}Let $\Gamma\Rightarrow \Delta$ be the end sequent of a $\bf{PLK}$-proof. Suppose also that it is the conclusion of an instance of some rule, say $R$. If $\Gamma\Rightarrow \Delta$ has a $\bf{PLK}$-proof, then each premise of $R$ also has a $\bf{PLK}$-proof.\end{thm}}

\hypertarget{ml2}{\begin{thm}{\label{ml2}}Let $\Gamma\Rightarrow \Delta$ be the end sequent of a $\bf{PLK}$-proof. Suppose also that it is the conclusion of an instance of a rule, say $R$. If this proof of $\Gamma\Rightarrow \Delta$ is \textsc{contraction+cut}-free, then each premise of $R$ also has a \textsc{contraction+cut}-free $\bf{PLK}$-proof.\end{thm}}

\hypertarget{thm:atomic}{\begin{thm}{\label{thm:atomic}}Every \textbf{PLK}-sequent of the form $A\Rightarrow  A$ has a \textsc{contraction+cut}-free proof.\end{thm}}

\begin{proof}
We use structural induction to prove this result. The base case is done by noting that $p\Rightarrow p$ is provable without using any instance of \textsc{contraction} or \textsc{cut} because the one-step proof $$p\Rightarrow p$$suffices.

Let us now assume that for $A,B\in \mathscr{L}_{\bf{CPL}}$, both $A\Rightarrow A$ and $B\Rightarrow B$ are provable. We now show that each of the sequents $\neg A\Rightarrow  \neg A, A\land B\Rightarrow  A\land B, A\lor B\Rightarrow  A\lor B$ and $A\to B\Rightarrow  A\to B$ is provable. 

For $\neg A\Rightarrow  \neg A$, after the proof of $A\Rightarrow A$, we do the following steps.

\begin{prooftree}
\def\fCenter{\Rightarrow}
\Axiom$A\fCenter A$
\RightLabel{{\color{red}{\quad $\neg\sf{L}$}}}
\UnaryInf$\neg A,A\fCenter \Lambda$
\RightLabel{{\color{red}{\quad $\sf{EL}$}}}
\UnaryInf$A,\neg A\fCenter \Lambda$
\RightLabel{{\color{red}{\quad $\neg\sf{R}$}}}
\UnaryInf$\neg A\fCenter \neg A$
\end{prooftree}

For $A\land B\Rightarrow  A\land B$, after the proofs of $A\Rightarrow A$ and $B\Rightarrow B$, we do the following steps.

\begin{prooftree}
\def\fCenter{\Rightarrow}
\Axiom$A\fCenter A$
\LeftLabel{{\color{red}{$\land\sf{L}$\quad}}}
\UnaryInf$A\land B\fCenter A$
\Axiom$B\fCenter B$
\RightLabel{{\color{red}{\quad $\land\sf{L}$}}}
\UnaryInf$A\land B\fCenter B$
\RightLabel{{\color{red}{\quad $\land\sf{R}$}}}
\BinaryInf$A\land B\fCenter A\land B$
\end{prooftree}

For $A\lor B\Rightarrow  A\lor B$, after the proofs of $A\Rightarrow A$ and $B\Rightarrow B$, we do the following steps.

\begin{prooftree}
\def\fCenter{\Rightarrow}
\Axiom$A\fCenter A$
\LeftLabel{{\color{red}{$\lor\sf{R}$\quad}}}
\UnaryInf$A\fCenter A\lor B$
\Axiom$B\fCenter B$
\RightLabel{{\color{red}{\quad $\lor\sf{R}$}}}
\UnaryInf$B\fCenter A\lor B$
\RightLabel{{\color{red}{\quad $\lor\sf{L}$}}}
\BinaryInf$A\lor B\fCenter A\lor B$
\end{prooftree}

For $A\to B\Rightarrow  A\to B$, after the proofs of $A\Rightarrow A$ and $B\Rightarrow B$, we do the following steps.

\begin{prooftree}
\def\fCenter{\Rightarrow}
\Axiom$A\fCenter A$
\LeftLabel{{\color{red}{$\sf{WR}$\quad}}}
\UnaryInf$A\fCenter A, B$
\LeftLabel{{\color{red}{$\sf{ER}$\quad}}}
\UnaryInf$A\fCenter B,A$
\Axiom$B\fCenter B$
\RightLabel{{\color{red}{\quad $\sf{WL}$}}}
\UnaryInf$A,B\fCenter B$
\RightLabel{{\color{red}{\quad $\sf{EL}$}}}
\UnaryInf$B,A\fCenter B$
\RightLabel{{\color{red}{\quad $\to\sf{L}$}}}
\BinaryInf$A\to B, A\fCenter B$
\RightLabel{{\color{red}{\quad $\sf{EL}$}}}
\UnaryInf$A,A\to B\fCenter B$
\RightLabel{{\color{red}{\quad $\to\sf{R}$}}}
\UnaryInf$A\to B\fCenter A\to B$
\end{prooftree}

Each of these proofs is \textsc{contraction+cut}-free. Hence, by the principle of structural induction, we are done.
\end{proof}

\hypertarget{thm:con+cut_elim(I)}{\begin{thm}[\textsc{Contraction-Elimination Theorem for \textbf{PLK}}]{\label{thm:con+cut_elim(I)}}Let $\Gamma\Rightarrow \Delta$ has a \textbf{PLK}-proof. Then, it has a \textsc{contraction}-free proof as well.\end{thm}}

\begin{proof}
Consider the highest leftmost application of \textsc{contraction}. Let the sequent obtained by this application be $\Gamma'\Rightarrow \Delta'$. Two cases arise.

\begin{case}
If $\Gamma'\Rightarrow \Delta'$ is the result of an application of {\color{red}{\sf{CL}}} then we have, 

\begin{prooftree}
\def\fCenter{\Rightarrow}
\Axiom$A, A,\Gamma''\fCenter \Delta'$
\RightLabel{{\color{red}{\quad $\sf{CL}$}}}
\UnaryInf$A,\Gamma''\fCenter \Delta'$
\end{prooftree}
where $A,\Gamma''\Rightarrow \Delta'$ is $\Gamma'\Rightarrow \Delta'$. Replace this with the following application of \textsc{cut}, 

\begin{prooftree}
\def\fCenter{\Rightarrow}
\AxiomC{$A,\Gamma''\fCenter \Delta', A$}
\AxiomC{$A, A,\Gamma''\fCenter \Delta'$}
\RightLabel{{\color{red}{\quad $\textsc{cut}$}}}
\BinaryInfC{$A,\Gamma''\fCenter \Delta'$}
\end{prooftree}
\end{case}

\begin{case}
If $\Gamma'\Rightarrow \Delta'$ is the result of an application of {\color{red}{\sf{CR}}} then we have, 

\begin{prooftree}
\def\fCenter{\Rightarrow}
\Axiom$\Gamma'\fCenter \Delta'',A,A$
\RightLabel{{\color{red}{\quad $\sf{CR}$}}}
\UnaryInf$\Gamma'\fCenter \Delta'', A$
\end{prooftree} where $\Gamma'\Rightarrow \Delta'', A$ is $\Gamma'\Rightarrow \Delta'', A$. Replace this by the following application of \textsc{cut}, 

\begin{prooftree}
\def\fCenter{\Rightarrow}
\AxiomC{$\Gamma'\fCenter \Delta'', A,A$}
\AxiomC{$A,\Gamma'\fCenter \Delta'', A$}
\RightLabel{{\color{red}{\quad $\textsc{cut}$}}}
\BinaryInfC{$\Gamma'\fCenter \Delta'', A$}
\end{prooftree}
\end{case}
The upper leftmost sequent (in the first case) and the upper rightmost sequent (in the second case) can be obtained by application(s) of \textsc{weakening} and/or \textsc{exchange} from $A\Rightarrow A$. Since the proof of $A\Rightarrow A$ can be done without using \textsc{contraction+cut} (by \hyperlink{thm:atomic}{\thmref{thm:atomic}}). Repeat this process until all the instances of \textsc{contraction} have been removed. 
\end{proof}

For the following theorem, we will need the notion of height of a proof. 
\begin{defn}[\textsc{Height of a \textbf{PLK}-proof}] The \textit{height of a $\mathbf{PLK}$-proof} $\mathcal{D}$, denoted as $h(\mathcal{D})$ is defined recursively as follows:
\begin{itemize}
    \item If $\mathcal{D}$ is the single-step proof of an axiom, then $h(\mathcal{D})=0$.
    \item If $\mathcal{D}$ is the proof of a sequent $S$, deduced from $S'$ by an application of a single-premise rule, with $\mathcal{D}'$ being a proof of $S'$, then $h(\mathcal{D})=h(\mathcal{D}')+1$.
    \item If $\mathcal{D}$ is the proof of a sequent $S$, deduced from $S'$ and $S''$ by an application of a two-premise rule, with $\mathcal{D}'$ and $\mathcal{D}''$ being their proofs respectively, then $h(\mathcal{D})=\max(h(\mathcal{D}',\mathcal{D}'')+1$.
\end{itemize}
\end{defn}

We will also need the notion of \textit{degree} of a proof. 

\begin{defn}[\textsc{Degree of a \textbf{PLK}-proof}]Given a $\bf{PLK}$ proof, the \textit{degree} of the proof is the maximum complexity (i.e., the number of connectives appearing in a formula) of the \textsc{cut}-formula appearing in the proof. 

A given sequent $\Gamma\Rightarrow \Delta$ is said to have a \textit{proof in $\mathbf{PLK}_n$}, more simply \textit{$\mathbf{PLK}_n$-proof}, if there is a $\bf{PLK}$-proof of $\Gamma\Rightarrow \Delta$, whose degree is at most $n$.\end{defn}

In \hyperlink{thm:PLK_cons}{\thmref{thm:PLK_cons}}, we will show that certain sequents are not \textbf{PLK}-provable. In particular, we will show that $\Lambda\Rightarrow \Lambda$ is not \textbf{PLK}-provable. Usually, this result is proved as a corollary of the \textsc{cut}-elimination theorem of \textbf{PLK} or via soundness theorem. In fact, both \hyperlink{thm:PLK0_cons}{\thmref{thm:PLK0_cons}} and \hyperlink{thm:PLK_cons}{\thmref{thm:PLK_cons}} admit easier proofs if one uses the soundness theorem for \textbf{PLK}. However, due to our reliance on purely syntactic methods, we can't apply the soundness theorem. Also, since we are trying to prove the \textsc{cut}-elimination theorem for \textbf{PLK}, we can't use it to prove either \hyperlink{thm:PLK0_cons}{\thmref{thm:PLK0_cons}} or \hyperlink{thm:PLK_cons}{\thmref{thm:PLK_cons}}.

\hypertarget{thm:PLK0_cons}{\begin{thm}{\label{thm:PLK0_cons}}There is no $\bf{PLK}_0$-proof of sequents of the following forms: 
\begin{enumerate}[label=(\arabic*)]
\item $\Lambda\Rightarrow \Lambda$
\item $p_1,\ldots,p_k\Rightarrow \Lambda$
\item $\Lambda\Rightarrow q_1,\ldots,q_l$
\item $p_1,\ldots,p_k\Rightarrow q_1,\ldots,q_l$
\end{enumerate} where $p_i$'s and $q_j$'s are propositional variables and for (4) $\{p_1,\ldots,p_k\}\cap \{q_1,\ldots,q_l\}=\emptyset$.\end{thm}}

\begin{proof}
Notice that, due to \hyperlink{thm:con+cut_elim(I)}{\thmref{thm:con+cut_elim(I)}}, it is enough to consider only \textsc{contraction}-free $\bf{PLK}_0$-proofs. We apply strong induction on the height of these proofs, i.e., on $h(\mathcal{D})$ where $\mathcal{D}$ is a \textsc{contraction}-free $\bf{PLK}_0$-proof.
    
$\\$\underline{\textsc{Base Case:}} In this case $h(\mathcal{D})=0$, and in this case the statement is true since the only proofs of height $0$ are of axiomatic sequents, which are of the form $r\Rightarrow r$, where $r$ is a propositional variable, and this is not among the forms listed in the theorem. 
    
$\\$\underline{\textsc{Induction Hypothesis:}} There is no \textsc{contraction}-free $\bf{PLK}_0$-proof of any sequent of the forms given in the theorem, of height $<n$.

$\\$\underline{\textsc{Induction Step:}} Consider any $\bf{PLK}_0$-proof of any one of the sequent of the form $\Lambda\Rightarrow \Lambda,p_1,\ldots,p_k\Rightarrow \Lambda, \Lambda\Rightarrow q_1,\ldots,q_l$ and $p_1,\ldots,p_k\Rightarrow q_1,\ldots,q_l$, where $p_i$'s and $q_j$'s are propositional variables, of height $n$ (and for the last case $\{p_1,\ldots,p_k\}\cap \{q_1,\ldots,q_l\}=\emptyset$). Four cases will arise, one for each sequent in the list. In each case, we will get a contradiction. We discuss them one by one. 

\begin{case}
    Suppose that the sequent under consideration is $\Lambda\Rightarrow \Lambda$ and that it has a \textbf{PLK}-proof of height $n$.

    Then, note that since the proof is \textsc{contraction}-free and since no logical connective is introduced in $\Lambda\Rightarrow \Lambda$, it can't be the conclusion of any instance of a logical rule. Furthermore, since no formula is introduced and no formulas are interchanged, it can't be the conclusion of any instance of \textsc{weakening} and \textsc{exchange}-rules. Hence $\Lambda\Rightarrow \Lambda$ can only be the conclusion of an instance of \textsc{cut}. Hence, the last couple of lines of the proof will look something like the following: 

    \begin{prooftree}
    \def\fCenter{\Rightarrow}
    \Axiom$\Lambda\fCenter p$
    \RightLabel{{\color{red}{\quad $\textsc{cut}$}}}
    \Axiom$p\fCenter \Lambda $
    \BinaryInf$\Lambda\fCenter \Lambda$
    \end{prooftree}
    where $p$ is a propositional variable (here we use the fact that this proof is a $\bf{PLK}_0$-proof).

    But, by \hyperlink{ml1}{\thmref{ml1}}, this would imply that there exists a $\bf{PLK}_0$-proof of $\Lambda\Rightarrow  p$ of height less than $n$, contrary to the induction hypothesis. Consequently, in this case, our conclusion follows.
\end{case}

\hypertarget{case:2(7)}{\begin{case}{\label{case:2(7)}}
    Suppose that the sequent under consideration is $p_1,\ldots,p_k\Rightarrow \Lambda$ and that it has a \textbf{PLK}-proof of height $n$.}

    Then, note that since the proof is \textsc{contraction}-free and since no logical connective is introduced in $p_1,\ldots,p_k\Rightarrow \Lambda$, it can only be the conclusion of an instance of \textsc{cut} or \textsc{weakening} or \textsc{exchange}. 

    Now, starting from the bottom, we move upward. Proceeding in this manner, at the first occurrence of a sequent which is not the conclusion of an instance of \textsc{exchange}. Let that sequent be $p_{\sigma(1)},\ldots,p_{\sigma(k)}\Rightarrow \Lambda$ where $\sigma:\{1,\ldots,k\}\to \{1,\ldots,k\}$ is a bijection. Then, this sequent should be the conclusion of an instance of \textsc{cut} or \textsc{weakening}.
    
    \begin{subcase} If $k>1$, then the following cases arise.
        \begin{subsubcase}        
         Suppose that $p_{\sigma(1)},\ldots,p_{\sigma(k)}\Rightarrow \Lambda$ is the conclusion of an instance of \textsc{weakening}. Then, it must be the conclusion of an instance of ${\color{red}{\mathsf{WL}}}$. Hence, the last two lines of the proof will look something like the following: 

        \begin{prooftree}
        \def\fCenter{\Rightarrow}
        \Axiom$p_{\sigma(2)},\ldots,p_{\sigma(k)}\fCenter \Lambda $
        \RightLabel{{\color{red}{\quad $\mathsf{WL}$}}}    \UnaryInf$p_{\sigma(1)},p_{\sigma(2)}\ldots,p_{\sigma(k)}\fCenter \Lambda $
        \end{prooftree}

        But then, this would mean that $p_{\sigma(2)},\ldots,p_{\sigma(k)}\Rightarrow \Lambda$ has a $\bf{PLK}_0$-proof of height $<n$, a contradiction to the induction hypothesis. Consequently, in this case also, our conclusion follows.
        \end{subsubcase}
    \begin{subsubcase}        
        If, on the other hand, it is the conclusion of an instance of \textsc{cut}, then the last two lines of the proof will look something like the following: 

        \begin{prooftree}
        \def\fCenter{\Rightarrow}        \Axiom$p_{\sigma(1)},p_{\sigma(2)},\ldots,p_{\sigma(k)}\fCenter r$ 
        \Axiom$r, p_{\sigma(1)},p_{\sigma(2)},\ldots,p_{\sigma(k)} \fCenter \Lambda$   
        \RightLabel{{\color{red}{\quad $\textsc{cut}$}}}    \BinaryInf$p_{\sigma(1)},p_{\sigma(2)}\ldots,p_{\sigma(k)}\fCenter \Lambda $
        \end{prooftree}
        For some propositional variable $r$. (Here again, we use the fact that this proof is a $\bf{PLK}_0$-proof).

        If $r\in \{p_{\sigma(1)},p_{\sigma(2)}\ldots,p_{\sigma(k)}\}$ then it would imply that the sequent $r, p_{\sigma(1)},\ldots,p_{\sigma(k)}\Rightarrow \Lambda$ has a $\bf{PLK}_0$-proof of height $<n$, contrary to the induction hypothesis. If $r\notin\{p_{\sigma(1)},p_{\sigma(2)}\ldots,p_{\sigma(k)}\}$ then by considering the sequent  $ p_{\sigma(1)},\ldots,p_{\sigma(k)}\Rightarrow r$ we will obtain a similar contradiction. Thus, we are done in this case as well.
        \end{subsubcase}
    \end{subcase}
    \begin{subcase}
        If $k=1$, one can also obtain a contradiction by slightly modifying the above argument. Hence, in this case also, we are done.
    \end{subcase}    
\end{case}

\begin{case}
    Suppose that the sequent under consideration is $\Lambda\Rightarrow q_1,\ldots,q_l$ and that it has a \textbf{PLK}-proof of height $n$. Then, this case can be dealt with analogously to the previous one.
\end{case}

\begin{case}
    Suppose that the sequent under consideration is $p_1,\ldots,p_k\Rightarrow q_1,\ldots,q_l$ such that $\{p_1,\ldots,p_k\}\cap\{q_1,\ldots,q_l\}=\emptyset$ and that it has a \textbf{PLK}-proof of height $n$. Then, in this case both $k,l\ge 1$. 

    Without loss of generality, we may assume that $p_1,\ldots,p_k\Rightarrow q_1,\ldots,q_l$ is not the conclusion of any instance of an \textsc{exchange}-rule (if it were, then we will use the argument as in \hyperlink{case:2(7)}{\caseref{case:2(7)}} to obtain a sequent which is not). 

    So, since $p_1,\ldots,p_k\Rightarrow q_1,\ldots,q_l$ is not the conclusion of any instance of an \textsc{exchange}-rule it must be the conclusion of some instance of \textsc{weakening} or \textsc{cut}. In the former case, the possibilities for the last two lines of the proof look something like the following:

    \begin{prooftree}
        \def\fCenter{\Rightarrow}        
        \Axiom$p_1,\ldots,p_k\fCenter q_1,\ldots,q_{l-1}$ 
        \RightLabel{{\color{red}{\quad $\mathsf{WR}$}}}   
        \UnaryInf$p_1,\ldots,p_k\fCenter q_1,\ldots,q_{l-1},q_l$   
    \end{prooftree}
    if $l>1$, or, 

    \begin{prooftree}
        \def\fCenter{\Rightarrow}        
        \Axiom$p_1,\ldots,p_k\fCenter \Lambda$ 
        \RightLabel{{\color{red}{\quad $\mathsf{WR}$}}}   
        \UnaryInf$p_1,\ldots,p_k\fCenter q_1$   
    \end{prooftree}
    if $l=1$, or,

    \begin{prooftree}
        \def\fCenter{\Rightarrow}        
        \Axiom$p_1,\ldots,p_k\fCenter q_1,\ldots,q_{l}$ 
        \RightLabel{{\color{red}{\quad $\mathsf{WL}$}}}   
        \UnaryInf$p_1,p_2,\ldots,p_k\fCenter q_1,\ldots,q_l$   
    \end{prooftree}
    if $k>1$, or, 

    \begin{prooftree}
        \def\fCenter{\Rightarrow}        
        \Axiom$\Lambda\fCenter q_1,\ldots,q_{l}$ 
        \RightLabel{{\color{red}{\quad $\mathsf{WL}$}}}   
        \UnaryInf$p_1\fCenter q_1,\ldots,q_{l}$   
    \end{prooftree}
    if $k=1$, each of which will imply that there exists a proof of height less than $n$ of a sequent of one of the types listed in (1), (2) and (3), a contradiction. This shows that $p_1,\ldots,p_k\Rightarrow q_1,\ldots,q_l$ can't be the conclusion of any instance of a \textsc{weakening}-rule. 
    If it is the conclusion of an instance of \textsc{cut}, then the last two lines of the proof will look something like the following.
    \begin{prooftree}
        \def\fCenter{\Rightarrow}        \Axiom$p_1,\ldots,p_k\fCenter q_1,\ldots,q_{l}, r$ 
        \Axiom$r, p_1,\ldots,p_k\fCenter  q_1,\ldots,q_{l}$   
        \RightLabel{{\color{red}{\quad $\textsc{cut}$}}}    \BinaryInf$p_1,\ldots,p_k\fCenter q_1,\ldots,q_{l}$
    \end{prooftree}
If $r\in \{p_1,\ldots,p_k\}$ then consider the sequent $r, p_1,\ldots,p_k\Rightarrow q_1,\ldots,q_l$ and note that it is among the forms listed, $\{r,p_1,\ldots,p_k\}\cap \{q_1,\ldots,q_l\}=\emptyset$ but have a proof of height $<n$, a contradiction. If, on the other hand $r\notin \{p_1,\ldots,p_k\}$, then consider the $p_1,\ldots,p_k\Rightarrow q_1,\ldots,q_l,r$ to obtain a contradiction. 
\end{case}
Since we have obtained a contradiction in each of the cases, this completes the induction step, and our conclusion follows.
\end{proof}

We now state the following theorem without proof. For a proof of this theorem see \cite[Theorem 2.2]{Indrzejczak2021}. 

\hypertarget{thm:invertibility}{\begin{thm}[\textsc{Syntactic Invertibility of \textbf{PLK}}]{\label{thm:invertibility}}In $\mathbf{PLK}$ the following results hold.
\begin{enumerate}[label=(\arabic*)]
    \hypertarget{thm:invertibility(1)}{\item{\label{thm:invertibility(1)}} If $\emptyset\vdash_{\textbf{PLK}_n}\Gamma\Rightarrow \Delta, \neg A$ then $\emptyset\vdash_{\textbf{PLK}_n}A,\Gamma\Rightarrow \Delta$.}
    \hypertarget{thm:invertibility(2)}{\item{\label{thm:invertibility(2)}} If $\emptyset\vdash_{\textbf{PLK}_n}\neg A,\Gamma\Rightarrow \Delta$ then $\emptyset\vdash_{\textbf{PLK}_n}\Gamma\Rightarrow \Delta, A$.}
    \hypertarget{thm:invertibility(3)}{\item{\label{thm:invertibility(3)}} If $\emptyset\vdash_{\textbf{PLK}_n}A\land B,\Gamma\Rightarrow \Delta$ then $\emptyset\vdash_{\textbf{PLK}_n}A,B,\Gamma\Rightarrow \Delta$.}
    \hypertarget{thm:invertibility(4)}{\item{\label{thm:invertibility(4)}} If $\emptyset\vdash_{\textbf{PLK}_n}\Gamma\Rightarrow \Delta, A\land B$ then $\emptyset\vdash_{\textbf{PLK}_n}\Gamma\Rightarrow \Delta, A$ as well as $\emptyset\vdash_{\textbf{PLK}_n}\Gamma\Rightarrow \Delta, B$.}
    \hypertarget{thm:invertibility(5)}{\item{\label{thm:invertibility(5)}} If $\emptyset\vdash_{\textbf{PLK}_n}A\lor B,\Gamma\Rightarrow \Delta$ then $\emptyset\vdash_{\textbf{PLK}_n}A,\Gamma\Rightarrow \Delta$ as well as $\emptyset\vdash_{\textbf{PLK}_n}B,\Gamma\Rightarrow \Delta$.}
    \hypertarget{thm:invertibility(6)}{\item{\label{thm:invertibility(6)}} If $\emptyset\vdash_{\textbf{PLK}_n}\Gamma\Rightarrow \Delta, A\lor B$ then $\emptyset\vdash_{\textbf{PLK}_n}\Gamma\Rightarrow \Delta, A, B$.}
    \hypertarget{thm:invertibility(7)}{\item{\label{thm:invertibility(7)}} If $\emptyset\vdash_{\textbf{PLK}_n}A\to B,\Gamma\Rightarrow \Delta$ then $\emptyset\vdash_{\textbf{PLK}_n}\Gamma\Rightarrow \Delta, A$ as well as $\emptyset\vdash_{\textbf{PLK}_n}B,\Gamma\Rightarrow \Delta$.}
    \hypertarget{thm:invertibility(8)}{\item{\label{thm:invertibility(8)}} If $\emptyset\vdash_{\textbf{PLK}_n}\Gamma\Rightarrow \Delta, A\to B$ then $\emptyset\vdash_{\textbf{PLK}_n}A,\Gamma\Rightarrow \Delta, B$.}
\end{enumerate}\end{thm}}

\begin{rem}
    This is not the original statement of \cite[Theorem 2.2]{Indrzejczak2021}. However, it follows immediately from the argument given there.  
\end{rem}

\hypertarget{rem:invertibility}{\begin{rem}\label{rem:invertibility}
    It is important to mention at this point that even though there are different ways to prove this particular theorem, we are interested only in the proofs using \textit{global} transformations (see, e.g. \cite{Curry1977} and \cite{NegriPlato2001}). These types of proofs have the advantage that if the proof of the former sequent doesn't have any instance of the application of a particular rule, the proof of the latter also will not.  
\end{rem}}

Using these two theorems, we can prove the following theorem, which says that \textbf{PLK} is consistent. 

\hypertarget{thm:PLK_cons}{\begin{thm}{\label{thm:PLK_cons}}There is no $\bf{PLK}$-proof of the sequents of the following forms: 
\begin{enumerate}[label=(\arabic*)]
\item $\Lambda\Rightarrow \Lambda$
\item $p_1,\ldots,p_k\Rightarrow \Lambda$
\item $\Lambda\Rightarrow q_1,\ldots,q_l$
\item $p_1,\ldots,p_k\Rightarrow q_1,\ldots,q_l$
\end{enumerate} where $p_i$'s and $q_j$'s are propositional variables and for (4) $\{p_1,\ldots,p_k\}\cap \{q_1,\ldots,q_l\}=\emptyset$.\end{thm}}

\begin{proof} We show this only for $\Lambda\Rightarrow \Lambda$. The proofs for the other cases are similar. 

We show that if there is a $\mathbf{PLK}_n$-proof of $\Lambda\Rightarrow\Lambda$, there is a $\mathbf{PLK}_0$-proof of the same. 

Start with the highest leftmost application of \textsc{cut}. Let $\Gamma\Rightarrow\Delta$ be the sequent obtained by this application of \textsc{cut}. Then, the last two lines of the this part of the proof, ending in $\Gamma\Rightarrow\Delta$ will look something like the following.
\begin{prooftree}
        \def\fCenter{\Rightarrow}        
        \Axiom$\Gamma\fCenter\Delta, A$ 
        \Axiom$A,\Gamma\fCenter \Delta$
        \RightLabel{{\color{red}{\quad $\textsc{cut}$}}}   
        \BinaryInf$\Gamma\fCenter \Delta$   
\end{prooftree}Then, it follows that there is a $\mathbf{PLK}_{\mathsf{Comp}(A)}$-proof of $\Gamma\Rightarrow\Delta$. Furthermore, by our choice of $\Gamma\Rightarrow\Delta$, both $\Gamma\Rightarrow\Delta, A$ and $A, \Gamma\Rightarrow\Delta$ has $\mathbf{PLK}_0$-proofs. If $\mathsf{Comp}(A)=0$, then we do nothing. Otherwise, we do the following. 

\begin{case}
    If $A=\neg B$ for some formula $B$. Then, by \hyperlink{thm:invertibility}{\thmref{thm:invertibility}}, we replace the proofs of $\Gamma\Rightarrow \Delta,\neg B$ and $\neg B,\Gamma\Rightarrow \Delta$ by the proofs of $ B, \Gamma\Rightarrow \Delta$ and $\Gamma\Rightarrow \Delta, B$ respectively. Then, we get the following proof of $\Gamma\Rightarrow \Delta$. 
    \begin{prooftree}
        \def\fCenter{\Rightarrow}        
        \Axiom$\Gamma\fCenter\Delta, B$ 
        \Axiom$B,\Gamma\fCenter \Delta$
        \RightLabel{{\color{red}{\quad $\textsc{cut}$}}}   
        \BinaryInf$\Gamma\fCenter \Delta$   
\end{prooftree}
\end{case}
\begin{case}
    Let $A=B\land C$ for some formulas $B$ and $C$. Then, by \hyperlink{thm:invertibility}{\thmref{thm:invertibility}}, we replace the proofs of $\Gamma\Rightarrow \Delta, B\land C$ and $B\land C,\Gamma\Rightarrow \Delta$ by the proofs of $ \Gamma\Rightarrow \Delta, B,\Gamma\Rightarrow \Delta, C$ and $B, C, \Gamma\Rightarrow \Delta$ respectively. Then, we get the following proof of $\Gamma\Rightarrow \Delta$. 
    \begin{prooftree}
        \def\fCenter{\Rightarrow} 
        \Axiom$\Gamma\fCenter\Delta, C$
        \Axiom$\Gamma\fCenter\Delta, B$ 
        \LeftLabel{{\color{red}{$\mathsf{WL}$\quad}}} 
        \UnaryInf$C,\Gamma\fCenter\Delta, B$
        \Axiom$B, C,\Gamma\fCenter\Delta$
        \RightLabel{{\color{red}{\quad $\textsc{cut}$}}}   
        \BinaryInf$C,\Gamma\fCenter\Delta$ 
        \RightLabel{{\color{red}{\quad $\textsc{cut}$}}}   
        \BinaryInf$\Gamma\fCenter\Delta$ 
\end{prooftree}
\end{case}

\begin{case}
    Let $A= B\lor C$ for some formulas $B$ and $C$. Then, by \hyperlink{thm:invertibility}{\thmref{thm:invertibility}}, we replace the proofs of $\Gamma\Rightarrow \Delta, B\lor C$ and $\ B\lor C, \Gamma\Rightarrow \Delta$ by the proofs of $B, \Gamma\Rightarrow \Delta, C,\Gamma\Rightarrow \Delta $ and $\Gamma\Rightarrow \Delta, B,C$ respectively. Then, we get the following proof of $\Gamma\Rightarrow \Delta$. 
    \begin{prooftree}
        \def\fCenter{\Rightarrow}         
        \Axiom$\Gamma\fCenter\Delta, B,C$
        \Axiom$C,\Gamma\fCenter\Delta$ 
        \RightLabel{{\color{red}{\quad $\mathsf{WR}$}}} 
        \UnaryInf$C,\Gamma\fCenter\Delta, B$
        \LeftLabel{{\color{red}{$\textsc{cut}$\quad }}}   
        \BinaryInf$\Gamma\fCenter\Delta, B$ 
        \Axiom$B,\Gamma\fCenter \Delta$
        \LeftLabel{{\color{red}{$\textsc{cut}$\quad }}} 
        \BinaryInf$\Gamma\fCenter \Delta$ 
\end{prooftree}
\end{case}

\begin{case}
    Let $A= B\to C$ for some formulas $B$ and $C$. Then, by \hyperlink{thm:invertibility}{\thmref{thm:invertibility}}, we replace the proofs of $\Gamma\Rightarrow \Delta, B\to C$ and $\ B\to C, \Gamma\Rightarrow \Delta$ by the proofs of $\Gamma\Rightarrow \Delta, B, C,\Gamma\Rightarrow\Delta$ and $B,\Gamma\Rightarrow\Delta, C$ respectively. Then, we get the following proof of $\Lambda\Rightarrow \Lambda$. 
    \begin{prooftree}
        \def\fCenter{\Rightarrow}         
        \Axiom$\Gamma\fCenter \Delta, B$
        \Axiom$B,\Gamma\fCenter \Delta, C$
        \Axiom$C, \Gamma\fCenter \Delta$ 
        \RightLabel{{\color{red}{\quad $\mathsf{WL}$}}} 
        \UnaryInf$B,C, \Gamma\fCenter \Delta$
        \RightLabel{{\color{red}{\quad $\mathsf{EL}$}}} 
        \UnaryInf$C,B, \Gamma\fCenter \Delta$
        \RightLabel{{\color{red}{\quad $\textsc{cut}$}}}
        \BinaryInf$B, \Gamma\fCenter \Delta$ 
        \RightLabel{{\color{red}{\quad $\textsc{cut}$}}}   
        \BinaryInf$\Gamma\fCenter \Delta$ 
\end{prooftree}
\end{case}
In each case, the complexity of the \textsc{cut}-formula is $<\mathsf{Comp}(A)$. Since $\mathsf{Comp}(A)\le n$, we thus obtain a proof of $\Gamma\Rightarrow\Delta$ of a lower degree than $n$. Continuing this process, till there is a \textsc{cut}-formula of nonzero degree, we will get a $\mathbf{PLK}_0$ proof of $\Lambda\Rightarrow\Lambda$, contradicting \hyperlink{thm:PLK0_cons}{\thmref{thm:PLK0_cons}}. This contradiction ensures that our assumption was wrong, and hence there is no $\mathbf{PLK}_n$-proof of $\Lambda\Rightarrow\Lambda$. Since $n$ was chosen arbitrarily, we have thus shown that there is no $\mathbf{PLK}$-proof of $\Lambda\Rightarrow\Lambda$.
\end{proof}

\hypertarget{ce1}{\begin{thm}[\textsc{Cut-Elimination Theorem for \textbf{PLK}}]{\label{ce1}}Let $\Gamma\Rightarrow \Delta$ be a $\mathbf{PLK}$-provable sequent. Then there is a \textsc{cut}-free $\bf{PLK}$-proof of $\Gamma\Rightarrow \Delta$. \end{thm}}

\begin{proof}
The proof is by induction on the value of $\SeqComp(\Gamma\Rightarrow \Delta)$, which denotes the total number of connectives in the wffs of $\Gamma$ and $\Delta$, $n$.

$\\$\underline{\textsc{Base Case:}} In the base case, $\SeqComp(\Gamma\Rightarrow \Delta) = 0$, the sequent $\Gamma\Rightarrow \Delta$ contains no logical connectives and thus every formula in the sequent is a propositional variable. Since the sequent is provable, by \hyperlink{thm:PLK_cons}{\thmref{thm:PLK_cons}}, there must be some variable $p$ which occurs both in $\Gamma$ and in $\Delta$. Thus $\Gamma\Rightarrow \Delta$ can be proved without using \textsc{cut}, although we may need \textsc{weakening} and \textsc{exchange} in this case.

$\\$\underline{\textsc{Induction Hypothesis:}} For every $\mathbf{PLK}$-sequent $\Gamma\Rightarrow \Delta$ if $\emptyset\vdash_{\textbf{PLK}} \Gamma\Rightarrow \Delta$ and $\SeqComp(\Gamma\Rightarrow \Delta)<m$, then there is a \textsc{cut}-free $\bf{PLK}$-proof of $\Gamma\Rightarrow \Delta$.

$\\$\underline{\textsc{Induction Step:}} Let $\Gamma\Rightarrow \Delta$ be a \textbf{PLK}-provable sequent for which $\SeqComp(\Gamma\Rightarrow \Delta)=m(>0)$. Without loss of generality (since we have \textsc{exchange} rule), we may assume that at least one of the extreme most formulas is not a propositional variable. Depending on the outermost connective of it, we have the following cases.

\begin{case}Let $\Gamma=\neg A,\Gamma'$ for some wff $A$. Then $\emptyset\vdash_{\textbf{PLK}}\neg A,\Gamma'\Rightarrow \Delta$. Hence, by \hyperlink{thm:invertibility}{\thmref{thm:invertibility}}  $\emptyset\vdash_{\textbf{PLK}}\Gamma'\Rightarrow \Delta, A$. Furthermore, since $\SeqComp(\Gamma'\Rightarrow \Delta, A)<m$ by the induction hypothesis, it follows that $\Gamma'\Rightarrow \Delta, A$ has a \textsc{cut}-free proof. Consequently, so does $\Gamma\Rightarrow \Delta$ (just apply ${\color{red}{\neg\,\sf{L}}}$ to $\Gamma'\Rightarrow \Delta, A$).
\end{case}

\begin{case}$\Delta=\Delta',\neg A$ for some wff $A$. This can be dealt in an analogous manner as the above.\end{case}

\begin{case}Let $\Gamma=A\lor B,\Gamma'$ for some wffs $A,B$. Then $\emptyset\vdash_{\textbf{PLK}}A\lor B,\Gamma'\Rightarrow \Delta$. Hence, by \hyperlink{thm:invertibility}{\thmref{thm:invertibility}}  $\emptyset\vdash_{\textbf{PLK}}A,\Gamma'\Rightarrow \Delta$ as well as $\emptyset\vdash_{\textbf{PLK}}B,\Gamma'\Rightarrow \Delta$. Furthermore, since $\SeqComp(A,\Gamma'\Rightarrow \Delta)<m$ as well as $\SeqComp(B,\Gamma'\Rightarrow \Delta)<m$ by the induction hypothesis, it follows that $A, \Gamma'\Rightarrow \Delta$ has a \textsc{cut}-free proof as well as $B, \Gamma'\Rightarrow \Delta$ has a \textsc{cut}-free proof. Consequently, so does $\Gamma\Rightarrow \Delta$ (just apply ${\color{red}{\lor\,\sf{R}}}$ to $A,\Gamma\Rightarrow \Delta$).\end{case}

\begin{case}Let $\Delta=\Delta', A\land B$ for some wffs $A,B$. This can be dealt in an analogous manner as the above.\end{case}

\begin{case}Let $\Gamma=A\to B,\Gamma'$ for some wffs $A,B$. This can be dealt in an analogous manner as the above.
\end{case}

\begin{case}Let $\Delta=\Delta', A\to B$ for some wffs $A,B$. This can be dealt in an analogous manner as the above.\end{case}

\begin{case}Let $\Delta=\Delta', A\lor B$ for some wffs $A,B$. Then  $\emptyset\vdash_{\textbf{PLK}}\Gamma\Rightarrow \Delta', A\lor B$. By \hyperlink{thm:invertibility}{\thmref{thm:invertibility}} $\emptyset\vdash_{\textbf{PLK}}\Gamma\Rightarrow \Delta', A, B$ as well. Furthermore, since $\SeqComp(\Gamma\Rightarrow \Delta', A, B)<m$ by the induction hypothesis, it follows that $\Gamma\Rightarrow \Delta', A, B$ has a \textsc{cut}-free proof. Consequently, $\Gamma\Rightarrow \Delta$ has a \textsc{cut}-free proof as follows:

\begin{prooftree}
\def\fCenter{\Rightarrow}
\Axiom$\Gamma\fCenter \Delta', A, B$
\RightLabel{{\color{red}{\quad $\lor\sf{R}$}}}
\UnaryInf$\Gamma\fCenter \Delta', A, A\lor B$
\RightLabel{{\color{red}{\quad $\sf{ER}$}}}
\UnaryInf$\Gamma\fCenter \Delta', A\lor B, A$
\RightLabel{{\color{red}{\quad $\lor\sf{R}$}}}
\UnaryInf$\Gamma\fCenter \Delta', A\lor B, A\lor B$
\RightLabel{{\color{red}{\quad $\sf{CR}$}}}
\UnaryInf$\Gamma\fCenter \Delta', A\lor B$
\end{prooftree}
Consequently, it follows that $\Gamma\Rightarrow \Delta$ has a \textsc{cut}-free proof.\end{case}

\begin{case}Let $\Gamma=A\land B,\Gamma'$ for some wffs $A,B$. Then  $\emptyset\vdash_{\textbf{PLK}}A\land B,\Gamma'\Rightarrow \Delta$. By \hyperlink{thm:invertibility}{\thmref{thm:invertibility}} $\emptyset\vdash_{\textbf{PLK}}A, B, \Gamma'\Rightarrow \Delta$ as well. Furthermore, $\SeqComp(A, B,\Gamma'\Rightarrow \Delta)<m$ by the induction hypothesis, it follows that $A, B,\Gamma\Rightarrow \Delta$ has a \textsc{cut}-free proof. $\Gamma\Rightarrow \Delta$ has a \textsc{cut}-free proof as follows:

\begin{prooftree}
\def\fCenter{\Rightarrow}
\Axiom$A, B,\Gamma'\fCenter \Delta$
\RightLabel{{\color{red}{\quad $\land\sf{L}$}}}
\UnaryInf$A\land B, B, \Gamma'\fCenter \Delta$
\RightLabel{{\color{red}{\quad $\sf{EL}$}}}
\UnaryInf$B, A\land B,\Gamma'\fCenter \Delta$
\RightLabel{{\color{red}{\quad $\land\sf{L}$}}}
\UnaryInf$A\land B, A\land B,\Gamma'\fCenter \Delta$
\RightLabel{{\color{red}{\quad $\sf{CL}$}}}
\UnaryInf$A\land B,\Gamma'\fCenter \Delta$
\end{prooftree}
Consequently, it follows that $\Gamma\Rightarrow \Delta$ has a \textsc{cut}-free proof.\end{case}
\end{proof}

\subsection{Second Proof of the \textsc{Cut}-elimination Theorem for \textbf{PLK}}
In this section, we provide a different argument to prove a slightly stronger result. We show that if a sequent has a \textbf{PLK}-proof, then it has a \textsc{contraction+cut}-free \textbf{PLK}-proof as well.  

\hypertarget{thm:prop_cont=>cont+cut}{\begin{thm}{\label{thm:prop_cont=>cont+cut}}Let $\Gamma\Rightarrow \Delta$ be a sequent and $q$ a propositional variable. If $q\to q,\Gamma\Rightarrow \Delta$ has a \textsc{contraction}-free proof in $\bf{PLK}$, then there is a \textsc{contraction+cut}-free $\bf{PLK}$-proof of $q\to q,\Gamma\Rightarrow \Delta$.\end{thm}}

\begin{proof}
We apply induction on the sequent complexity of $q\to q,\Gamma\Rightarrow \Delta$, $\SeqComp(q\to q,\Gamma\Rightarrow \Delta)$. 

$\\$\underline{\textsc{Base Case:}} In the base case, $\SeqComp(q\to q, \Gamma\Rightarrow \Delta) = 1$.  Since the sequent $q\to q,\Gamma\Rightarrow \Delta$ contains exactly one logical connective, and since $\to$ is a connective already occurring, it must be that one. Furthermore, every other formula in the sequent must be a propositional variable. 

Now, note that since $q\to q, \Gamma\Rightarrow \Delta$ is provable without using \textsc{contraction}, so are both the sequents $q, \Gamma\Rightarrow \Delta$ and $ \Gamma\Rightarrow \Delta, q$ (by \hyperlink{thm:invertibility}{\thmref{thm:invertibility}} and \hyperlink{rem:invertibility}{\remref{rem:invertibility}}). We claim that there exists a propositional variable $r$ occurring in both $\Gamma$ and $\Delta$. 

To show this, observe that since $q, \Gamma\Rightarrow \Delta$ is provable, by \hyperlink{thm:PLK_cons}{\thmref{thm:PLK_cons}}, it follows that there exists a propositional variable $r_1$ belonging to both the lists $q,\Gamma$ and $\Delta$. Similarly, since $\Gamma\Rightarrow \Delta,q$ is provable, by \hyperlink{thm:PLK_cons}{\thmref{thm:PLK_cons}}, it follows that there exists a propositional variable $r_2$ belonging to both the lists $\Gamma$ and $\Delta,q$. 

If $r_1\ne q$ or $r_2\ne q$, then by applying \textsc{weakening} and/or \textsc{exchange} as required, we obtain \textsc{contraction+cut}-free proofs of $q, \Gamma\Rightarrow \Delta$ and $ \Gamma\Rightarrow \Delta, q$. If, on the other hand, $r_1=r_2=q$, then since $\Delta$ contains $r_1$, it follows that $\Delta$ contains $q$; similarly, since $\Gamma$ contains $r_2$, it follows that $\Gamma$ contains $q$. Hence, in this case also, by applying \textsc{weakening} and/or \textsc{exchange} as required, we can obtain \textsc{contraction+cut}-free proofs of $q, \Gamma\Rightarrow \Delta$ and $ \Gamma\Rightarrow \Delta, q$.

The sequent $q\to q, \Gamma\Rightarrow \Delta$ can now be obtained by an application of $\color{red}{\to\mathsf{L}}$. We have thus obtained a \textsc{contraction+cut}-free proof of $q\to q, \Gamma\Rightarrow \Delta$. This completes the proof of the base case. 

$\\$\underline{\textsc{Induction Hypothesis:}} Suppose that for every sequent $\Gamma\Rightarrow \Delta$ with $\SeqComp(q\to q, \Gamma\Rightarrow \Delta)<n$ if there is a \textsc{contraction}-free \textbf{PLK}-proof of $q\to q, \Gamma\Rightarrow \Delta$, then there is a \textsc{contraction+cut}-free proof of $q\to q, \Gamma\Rightarrow \Delta$ as well.

$\\$\underline{\textsc{Induction Step:}} Consider a sequent $\Gamma\Rightarrow \Delta$ with $\SeqComp(q\to q, \Gamma\Rightarrow \Delta)=n$ such that $q\to q, \Gamma\Rightarrow \Delta$ is \textbf{PLK}-provable without using \textsc{contraction}. Consequently, by \hyperlink{thm:invertibility}{\thmref{thm:invertibility}} and \hyperlink{rem:invertibility}{\remref{rem:invertibility}}, it follows that both the sequents $q, \Gamma\Rightarrow \Delta$ and $\Gamma\Rightarrow \Delta,q$ are \textbf{PLK}-provable without using \textsc{contraction}. 

Depending upon the outermost connective occurring in the sequent $\Gamma\Rightarrow \Delta$, there will be eight cases -- four for the right and four for the left. Below, we discuss only the case when $\Gamma=B\lor C,\Gamma'$ for some formulas $B$ and $C$. Proofs for the other cases are similar.

Since $\Gamma=B\lor C,\Gamma'$, and since both $q, \Gamma\Rightarrow \Delta$ and $\Gamma\Rightarrow \Delta,q$ are \textbf{PLK}-provable without using \textsc{contraction}, it follows that both $q, B\lor C,\Gamma'\Rightarrow \Delta$ and $B\lor C,\Gamma'\Rightarrow \Delta,q$ are \textbf{PLK}-provable without using \textsc{contraction}. Consider $q, B\lor C,\Gamma'\Rightarrow \Delta$ and note that since it is provable, so is $B\lor C, q,\Gamma'\Rightarrow \Delta$ (apply $\color{red}{\mathsf{EL}}$). Then, by \hyperlink{thm:invertibility}{\thmref{thm:invertibility}} and \hyperlink{rem:invertibility}{\remref{rem:invertibility}}, it follows that there are \textsc{contraction}-free proofs of $B, q,\Gamma'\Rightarrow \Delta$ and $C, q,\Gamma'\Rightarrow \Delta$. Similarly, from the \textsc{contraction}-free proof of $B\lor C,\Gamma'\Rightarrow \Delta,q$, we will be able to construct \textsc{contraction}-free proofs of $B,\Gamma'\Rightarrow \Delta,q$ and $C,\Gamma'\Rightarrow \Delta,q$. Now, choose any such proof for each sequent. We extend each of these proofs as follows.
\begin{prooftree}
\def\fCenter{\Rightarrow}
\Axiom$B,\Gamma'\fCenter \Delta,q$
\Axiom$B,q,\Gamma'\fCenter \Delta$
\RightLabel{{\color{red}{\quad $\sf{EL}$}}}
\UnaryInf$q,B,\Gamma'\fCenter \Delta$
\RightLabel{{\color{red}{\quad$\to\sf{L}$}}}
\BinaryInf$q\to q, B,\Gamma'\fCenter \Delta$
\end{prooftree}
and,
\begin{prooftree}
\def\fCenter{\Rightarrow}
\Axiom$C,\Gamma'\fCenter \Delta,q$
\Axiom$C,q,\Gamma'\fCenter \Delta$
\RightLabel{{\color{red}{\quad $\sf{EL}$}}}
\UnaryInf$q,C,\Gamma'\fCenter \Delta$
\RightLabel{{\color{red}{\quad $\to\sf{L}$}}}
\BinaryInf$q\to q, C,\Gamma'\fCenter \Delta$
\end{prooftree}
Now, since $\SeqComp(q\to q, B,\Gamma'\Rightarrow \Delta)<n$ as well as $\SeqComp(q\to q, C,\Gamma'\Rightarrow \Delta)<n$ and since both are provable in \textbf{PLK} without using \textsc{contraction}, by the induction hypothesis it follows that each of these sequents has a \textsc{contraction+cut}-free proof in \textbf{PLK}. Choose any such for each of them. To these proofs, we add the following lines.

\begin{prooftree}
\def\fCenter{\Rightarrow}
\Axiom$q\to q, B,\Gamma'\fCenter \Delta$
\RightLabel{{\color{red}{\quad $\sf{EL}$}}}
\UnaryInf$B,q\to q,\Gamma'\fCenter \Delta$
\Axiom$q\to q, C,\Gamma'\fCenter \Delta$
\RightLabel{{\color{red}{\quad $\sf{EL}$}}}
\UnaryInf$C,q\to q,\Gamma'\fCenter \Delta$
\RightLabel{{\color{red}{\quad $\lor\sf{L}$}}}
\BinaryInf$B\lor C, q\to q,\Gamma'\fCenter \Delta$
\RightLabel{{\color{red}{\quad $\sf{EL}$}}}
\UnaryInf$q\to q, B\lor C,\Gamma'\fCenter \Delta$
\UnaryInf$q\to q, \Gamma\fCenter \Delta$
\end{prooftree}
Which gives a \textsc{contraction+cut}-free proof of $q\to q,\Gamma\Rightarrow \Delta$. This completes the proof in this case.
\end{proof}

\hypertarget{thm:con+cut_elim(II)}{\begin{thm}[\textsc{Contraction+Cut-Elimination Theorem for \textbf{PLK} (Part I)}]{\label{thm:con+cut_elim(II)}}Let $\Gamma\Rightarrow \Delta$ be a sequent and $A$ be a formula. If $A\to A,\Gamma\Rightarrow \Delta$ has a \textsc{contraction}-free proof in $\bf{PLK}$, then there is a \textsc{contraction+cut}-free $\bf{PLK}$-proof of $A\to A,\Gamma\Rightarrow \Delta$.\end{thm}}

\begin{proof}
The proof is done by induction on the complexity of $A$, $\mathsf{Comp}(A)$. The base case, i.e. when $\mathsf{Comp}(A)=0$, is proved in the previous theorem. The proof of the induction step follows an argument similar to the induction step of the last theorem.
\end{proof}

A \textit{multiset} may be formally defined as an ordered pair $(A, m)$ where $A$ is a set, called the \textit{underlying set} of the multiset and $m:A\to \mathbb{N}$, called the \textit{multiplicity function}. For any $a\in A$, $m(a) = \text{the multiplicity of}~a$, is the number of occurrences of $a$ in the multiset $A$ (e.g., if $a\notin A$, we have $m(a)=0$). 

Given two multisets $(A,m)$ and $(B,n)$, we will say that $(A,m)$ is \textit{included in} $(B,n)$ if $A\subseteq B$ and for all $a\in A$ $m(a)\le n(a)$. In this case, we will write $(A,m)\sqsubseteq (B,n)$. We will write $(A,m)\sqsubset (B,n)$ if $(A,m)\sqsubseteq (B,n)$ and there exists $a\in A$ such that $m(a)<n(a)$.

To each ordered list, one can associate a multiset by simply ignoring the order. Given a sequent $\Gamma\Rightarrow \Delta$, we will denote the multisets corresponding to $\Gamma$ and $\Delta$ as $\MultSet(\Gamma)$ and $\MultSet(\Delta)$ respectively.  

\begin{defn}[\textsc{Subsequent}]
Let $\Gamma\Rightarrow \Delta$ be a sequent. A sequent $\Gamma'\Rightarrow \Delta'$ is said to be a \textit{subsequent} of $\Gamma\Rightarrow \Delta$ if $\MultSet(\Gamma')\sqsubseteq \MultSet(\Gamma)$ and $\MultSet(\Delta')\sqsubseteq \MultSet(\Delta)$. In this case, we will write $\Gamma'\Rightarrow \Delta')\preceq \Gamma\Rightarrow \Delta)$. We will write $\Gamma'\Rightarrow \Delta'\prec \Gamma\Rightarrow \Delta$ if either $\MultSet(\Gamma')\sqsubset \MultSet(\Gamma)$ or $\MultSet(\Delta')\sqsubset \MultSet(\Delta)$.
\end{defn}

\hypertarget{subsequent}{\begin{thm}{\label{subsequent}}
Let $\Gamma\Rightarrow \Delta$ be a sequent. If it is deducible from a sequent $\Gamma'\Rightarrow \Delta'$ by application(s) of either \textsc{weakening} or \textsc{exchange}, then $\Gamma'\Rightarrow \Delta'$ is a subsequent of $\Gamma\Rightarrow \Delta$. 
\end{thm}}

\begin{proof}
Suppose that $\Gamma'\Rightarrow \Delta'$ is a sequent $\Gamma\Rightarrow \Delta$ is deducible from $\Gamma'\Rightarrow \Delta'$ by an appropriate number of applications of \textsc{weakening} and/or \textsc{exchange}. Let the part of such a  proof from $\Gamma'\Rightarrow \Delta'$ to $\Gamma\Rightarrow \Delta$ look something like the following.
\begin{prooftree}
\def\fCenter{\Rightarrow}
\Axiom$\Gamma_0\fCenter \Delta_0$
\UnaryInf$\Gamma_1\fCenter \Delta_1$
\end{prooftree}
$$\vdots$$
\begin{prooftree}
\def\fCenter{\Rightarrow}
\AxiomC{{\color{white}{$\Gamma\fCenter \Delta$}}}
\UnaryInf$\Gamma_k\fCenter \Delta_k$
\end{prooftree}
Where $\Gamma_0=\Gamma',\Gamma_k=\Gamma,\Delta_0=\Delta'$ and $\Delta_k=\Delta$ and for each $i\in \{0,\ldots, k\}$ $\Gamma_{i+1}\Rightarrow \Delta_{i+1}$ is obtained from $\Gamma_{i}\Rightarrow \Delta_{i}$ either by \textsc{weakening} or \textsc{exchange}.

Then note that for all $i\in \{0,\ldots, k\}$ if  $\Gamma_{i+1}\Rightarrow \Delta_{i+1}$ is obtained via an application of \textsc{weakening} from $\Gamma_{i}\Rightarrow \Delta_{i}$, we have $\Gamma_{i}\Rightarrow \Delta_{i}\prec\Gamma_{i+1}\Rightarrow \Delta_{i+1}$. Else, if it is obtained via an application of \textsc{exchange}, we have $\Gamma_{i}\Rightarrow \Delta_{i}=\Gamma_{i+1}\Rightarrow \Delta_{i+1}$. In either case, at each step, we have $\Gamma_{i}\Rightarrow \Delta_{i}\preceq\Gamma_{i+1}\Rightarrow \Delta_{i+1}$, which proves our claim.
\end{proof}

\hypertarget{adm_neg_free}{\begin{thm}{\label{adm_neg_free}}
Let $\Gamma\Rightarrow \Delta$ and $\Gamma'\Rightarrow \Delta'$ be two \textbf{PLK}-sequents. If there is a derivation of $A\to A, \Gamma\Rightarrow \Delta$ from $A\to A, \Gamma'\Rightarrow \Delta'$ using only \textsc{weakening} and/or \textsc{exchange}, then there is a derivation of $\Gamma\Rightarrow \Delta, A$ from $\Gamma'\Rightarrow \Delta', A$ using only \textsc{weakening} and/or \textsc{exchange}.
\end{thm}}

\begin{proof}
Immediate by applying induction on the total number of wffs in $\Gamma$ and $\Delta$.
\end{proof}

\hypertarget{thm:con+cut_elim(III)}{\begin{thm}[\textsc{Contraction+Cut-Elimination Theorem for \textbf{PLK} (Part II)}]{\label{thm:con+cut_elim(III)}}Let $\Gamma\Rightarrow \Delta$ be a sequent. If $A\to A,\Gamma\Rightarrow \Delta$ has a \textsc{contraction+cut}-free $\bf{PLK}$-proof, then so does $\Gamma\Rightarrow \Delta$. \end{thm}}

\begin{proof}
The proof is by induction on the values of $\SeqComp(A\to A,\Gamma\Rightarrow \Delta)$.

$\\$\underline{\textsc{Base Case:}} In the base case, $\SeqComp(A\to A, \Gamma\Rightarrow \Delta) = 1$, the sequent $A\to A,\Gamma\Rightarrow \Delta$ contains exactly one logical connective and since $\to$ is a connective already occurring, it must be that one. Hence, $A\to A$ is actually of the form $q\to q$, where $q$ is a propositional variable. Also, every other formula in the sequent must be a propositional variable. 

Furthermore, observe that in this proof, there is no application of \textsc{cut} as well as \textsc{contraction}. Hence, the only potentially applicable rules are the logical rules, \textsc{weakening} and \textsc{exchange}. However, the application of any logical rule other than $\color{red}{\to\sf{L}}$ is prohibited because that will increase sequent complexity, and \textsc{exchange} doesn't introduce any new formula. So, if we consider the place where $q\to q$ appeared first in the proof, we are forced to conclude that this must be due to the application of either an instance of \textsc{weakening} or of $\color{red}{\to\sf{L}}$.

In the former case, the sequent just above $A\to A, \Gamma\Rightarrow\Delta$ is $\Gamma\Rightarrow \Delta$. Since $A\to A, \Gamma\Rightarrow\Delta$ is \textbf{PLK}-provable, so is $ \Gamma\Rightarrow\Delta$. Hence, there exists some propositional variable $r$ common to both $\Gamma$ and $\Delta$ (by \hyperlink{thm:PLK_cons}{\thmref{thm:PLK_cons}}). So, by applying \textsc{weakening} and/or \textsc{exchange}, we can prove $\Gamma\Rightarrow \Delta$. 

In the latter case, following an argument similar to the base case of \hyperlink{thm:prop_cont=>cont+cut}{\thmref{thm:prop_cont=>cont+cut}}, we can prove that there is a propositional variable common to both $\Gamma$ and $\Delta$. Then, applying \textsc{weakening} and/or \textsc{exchange}, we can prove $\Gamma\Rightarrow \Delta$ in this case as well. Consequently, the base case is proved.

$\\$\underline{\textsc{Induction Hypothesis:}} For every sequent $\Gamma\Rightarrow \Delta$ for which $\SeqComp(A\to A,\Gamma\Rightarrow \Delta)<m$, if there is a \textsc{contraction+cut}-free $\bf{PLK}$-proof of $A\to A,\Gamma\Rightarrow \Delta$, then there is a \textsc{contraction+cut}-free $\bf{PLK}$-proof of $\Gamma\Rightarrow \Delta$.

$\\$\underline{\textsc{Induction Step:}} Let $\Gamma\Rightarrow \Delta$ be a sequent for which $\SeqComp(A\to A,\Gamma\Rightarrow \Delta)=m$. Suppose also that $A\to A, \Gamma\Rightarrow \Delta$  has a \textsc{contraction+cut}-free $\bf{PLK}$-proof. 

We do the following to transform this particular proof of $A\to A, \Gamma\Rightarrow \Delta$ to obtain a \textsc{contraction+cut}-free proof of $\Gamma\Rightarrow \Delta$. Starting from the bottom, we go upwards and find the sequent where a logical rule is used to obtain the conclusion sequent, say $\Gamma'\Rightarrow \Delta'$. This means, in particular, that every sequent below $\Gamma'\Rightarrow \Delta'$, in the proof of $A\to A, \Gamma\Rightarrow \Delta$, is obtained either by \textsc{weakening} or \textsc{exchange}. So, by \hyperlink{subsequent}{\thmref{subsequent}}, it follows that $\Gamma'\Rightarrow \Delta'\preceq A\to A, \Gamma\Rightarrow \Delta$.

We divide the argument into the following cases.
\begin{enumerate}[label=\textbf{\arabic*.}]
    \item The first appearance of $A\to A$ in $A\to A, \Gamma\Rightarrow \Delta$ is after $\Gamma'\Rightarrow \Delta'$. (\hyperlink{case:1}{\caseref{case:1}})
    \item The first appearance of $A\to A$ in $A\to A, \Gamma\Rightarrow \Delta$ is in or before $\Gamma'\Rightarrow \Delta'$. (\hyperlink{case:2}{\caseref{case:2}}) 
    \begin{enumerate}[label=\textbf{2.\arabic*.}]
        \item $A\to A$ occurs in $\Gamma'\Rightarrow \Delta'$. In this case, let $\Gamma^\ast$ be such that $\Gamma':=A\to A,\Gamma^\ast$. (\hyperlink{subcase:2.1}{\subcaseref{subcase:2.1}})  
        \begin{enumerate}[label=\textbf{2.1.\arabic*.}]
            \item $\SeqComp(A\to A, \Gamma^\ast\!\Rightarrow \Delta')=\SeqComp(A\to A, \Gamma\!\Rightarrow \Delta)$. (\hyperlink{subsubcase:2.1.1}{\subsubcaseref{subsubcase:2.1.1}}) 
            \item $\SeqComp(A\to A, \Gamma^\ast\!\Rightarrow \Delta')<\SeqComp(A\to A, \Gamma\!\Rightarrow \Delta)$. (\hyperlink{subsubcase:2.1.2}{\subsubcaseref{subsubcase:2.1.2}}) 
        \end{enumerate}
        \item $A\to A$ occurs before $\Gamma'\Rightarrow \Delta'$. (\hyperlink{subcase:2.2}{\subcaseref{subcase:2.2}})
            \begin{enumerate}[label=\textbf{2.2.\arabic*.}]
                \item $A\to A$ occurs in $\Gamma'\Rightarrow \Delta'$. (\hyperlink{subsubcase:2.2.1}{\subsubcaseref{subsubcase:2.2.1}})
                \item $A\to A$ does not occur in $\Gamma'\Rightarrow \Delta'$. (\hyperlink{subsubcase:2.2.2}{\subsubcaseref{subsubcase:2.2.2}})
            \end{enumerate}
        \end{enumerate}
    \end{enumerate}

\hypertarget{case:1}{\begin{case}{\label{case:1}}}If the first appearance of $A\to A$ in $A\to A, \Gamma\Rightarrow \Delta$ is after $\Gamma'\Rightarrow \Delta'$, then it must be the result of an application of a \textsc{weakening}-rule (since \textsc{exchange} doesn't introduce new formulas). Then, that part of the proof looks something like the following.

\begin{prooftree}
\def\fCenter{\Rightarrow}
\Axiom$\Gamma''\fCenter \Delta''$
\RightLabel{{\color{red}{\quad $\sf{WL}$}}}
\UnaryInf$A\to A, \Gamma''\fCenter\Delta''$
\end{prooftree}
For some sequent $\Gamma''\Rightarrow \Delta''$ occurring in the proof of $A\to A, \Gamma\Rightarrow \Delta$. Consequently, we have a \textsc{contraction+cut}-free $\bf{PLK}$-proof of $\Gamma''\Rightarrow \Delta''$. 

Finally, we apply \textsc{weakening} and \textsc{exchange} rules, as required, to obtain a \textsc{contraction+cut}-free $\bf{PLK}$-proof of $\Gamma\Rightarrow \Delta$
\end{case}

\hypertarget{case:2}{\begin{case}{\label{case:2}}}So suppose that the appearance of $A\to A$ in $A\to A, \Gamma\Rightarrow \Delta$ is in or before $\Gamma'\Rightarrow \Delta'$. 

\hypertarget{subcase:2.1}{\begin{subcase}{\label{subcase:2.1}}}
If $A\to A$ occurs in $\Gamma'\Rightarrow \Delta'$, then this is possible only if $\Gamma'\Rightarrow \Delta'$ has been obtained by an application of ${\color{red}{\to\sf{L}}}$. Then $\Gamma'=A\to A, \Gamma^\ast$. In this case, the line where it is applied will look something like the following.
\begin{prooftree}
\def\fCenter{\Rightarrow}
\AxiomC{$\Gamma^\ast\fCenter \Delta', A$}
\AxiomC{$A,\Gamma^\ast\fCenter \Delta'$}
\RightLabel{{\color{red}{\quad $\to\,\sf{L}$}}}
\BinaryInfC{$A\to A, \Gamma^\ast\fCenter\Delta'$}
\end{prooftree}
For some sequent $\Gamma^\ast\Rightarrow \Delta'$ occurring in the proof of $A\to A, \Gamma\Rightarrow \Delta$. Then, two cases arise.

\hypertarget{subsubcase:2.1.1}{\begin{subsubcase}{\label{subsubcase:2.1.1}}}
$\SeqComp(A\to A,\Gamma^\ast\Rightarrow \Delta')=\SeqComp(A\to A,\Gamma\Rightarrow \Delta)$. The difference between the sequents $A\to A,\Gamma^\ast\Rightarrow \Delta'$ and $A\to A,\Gamma\Rightarrow \Delta$ lies in the number of propositional variables only). Hence, $\Gamma^\ast\Rightarrow \Delta',A\preceq \Gamma\Rightarrow \Delta, A$. Consequently, by \hyperlink{adm_neg_free}{\thmref{adm_neg_free}}, we can easily obtain a proof of  $\Gamma\Rightarrow \Delta, A$ by appropriate applications of \textsc{weakening} and/or \textsc{exchange}.
\end{subsubcase}
\hypertarget{subsubcase:2.1.2}{\begin{subsubcase}{\label{subsubcase:2.1.2}}}
$\SeqComp(A\to A,\Gamma^\ast\Rightarrow \Delta')<\SeqComp(A\to A,\Gamma\Rightarrow \Delta)$. Then, since $A\to A, \Gamma^\ast\Rightarrow \Delta'$ has a \textsc{contraction+cut}-free proof, by the induction hypothesis, we can conclude that there exists a \textsc{contraction+cut}-free proof of $\Gamma^\ast\Rightarrow \Delta'$ as well. Now, observe that, \begin{align*}A\to A, \Gamma^\ast\Rightarrow \Delta'\preceq A\to A, \Gamma\Rightarrow \Delta&\implies \Gamma^\ast\Rightarrow \Delta'\preceq \Gamma\Rightarrow \Delta\end{align*}Consequently, since $\Gamma^\ast\Rightarrow \Delta'$ has a \textsc{contraction+cut}-free proof and since by \hyperlink{adm_neg_free}{\thmref{adm_neg_free}} $\Gamma\Rightarrow \Delta$ can be obtained from $\Gamma^\ast\Rightarrow \Delta'$ by applications of \textsc{weakening} and/or \textsc{exchange}; it follows that $\Gamma\Rightarrow \Delta$ has a \textsc{contraction+cut}-free proof as well.
\end{subsubcase}
\end{subcase}

\hypertarget{subcase:2.2}{\begin{subcase}{\label{subcase:2.2}}}
If $A\to A$ occurs before $\Gamma'\Rightarrow \Delta'$, then two cases arise.
\hypertarget{subsubcase:2.2.1}{\begin{subsubcase}{\label{subsubcase:2.2.1}}}
Suppose that $A\to A$ occurs in $\Gamma'\Rightarrow \Delta'$ as well. Suppose further that it is not due to any application of ${\color{red}{\to\sf{L}}}$ (because otherwise by an argument similar to the one in \hyperlink{subcase:2.1}{\subcaseref{subcase:2.1}}, we are done). Then some other rule must have been applied in the previous step. Since it must be a logical rule, the following seven cases arise, one for each rule .
    \begin{subsubsubcase}    
    $\Gamma'\Rightarrow \Delta'$ is the result of an application of ${\color{red}{\land\sf{L}}}$. Then that part of the proof looks something like the following.
    \begin{prooftree}
    \def\fCenter{\Rightarrow}
    \Axiom$B,\Gamma^\ast\fCenter \Delta'$
    \RightLabel{{\color{red}{\quad $\land\sf{L}$}}}
    \UnaryInf$B\land C, \Gamma^\ast\fCenter\Delta'$
    \end{prooftree}Where $\Gamma'=B\land C,\Gamma^\ast$. Also note that since $A\to A$ was introduced earlier in the proof, it must be in the list $\Gamma^\ast$ (otherwise, refer to \hyperlink{subsubcase:2.2.2}{\subsubcaseref{subsubcase:2.2.2}}). Let $\Gamma^\ast:=\Gamma_1,A\to A, \Gamma_2$. Hence $\Gamma'=B\land C,\Gamma^\ast=B\land C,\Gamma_1,A\to A, \Gamma_2$. Now, \begin{align*}&{\color{white}{\implies\,\,\,}}\Gamma'\Rightarrow \Delta'\preceq A\to A,\Gamma\Rightarrow \Delta\\&\implies B\land C,\Gamma^\ast\Rightarrow \Delta'\preceq A\to A,\Gamma\Rightarrow \Delta\\&\implies B\land C,\Gamma_1,A\to A, \Gamma_2\Rightarrow \Delta'\preceq A\to A,\Gamma\Rightarrow \Delta\\&\implies A\to A, B\land C,\Gamma_1, \Gamma_2\Rightarrow \Delta'\preceq A\to A,\Gamma\Rightarrow \Delta\\&\implies B\land C,\Gamma_1, \Gamma_2\Rightarrow \Delta'\preceq \Gamma\Rightarrow \Delta\end{align*}In other words,  $B\land C,\Gamma_1,\Gamma_2\Rightarrow \Delta'$ is a subsequent of $\Gamma\Rightarrow \Delta$.

    Now note that since $B\land C, \Gamma^\ast\Rightarrow \Delta'$ has a \textsc{contraction+cut}-free proof, so does $B, \Gamma^\ast\Rightarrow \Delta'$ (by \hyperlink{thm:invertibility}{\thmref{thm:invertibility}}). Now we extend this proof of $B, \Gamma^\ast\Rightarrow \Delta'$ in the following way.
    \begin{prooftree}
    \def\fCenter{\Rightarrow}
    \Axiom$B,\Gamma^\ast\fCenter \Delta'$
    \UnaryInf$B,\Gamma_1,A\to A,\Gamma_2\fCenter\Delta'$
    \RightLabel{{\color{red}{\quad $\sf{EL}$}}}
    \UnaryInf$A\to A,B,\Gamma_1,\Gamma_2\fCenter\Delta'$
    \end{prooftree}
    This, now constitutes a \textsc{contraction+cut}-free proof of $A\to A,B,\Gamma_1,\Gamma_2\Rightarrow \Delta'$. Furthermore, since we have, $\SeqComp(A\to A, B, \Gamma_1,\Gamma_2\Rightarrow \Delta')<\SeqComp(B\land C,\Gamma_1,A\to A,\Gamma_2\Rightarrow \Delta')\le \SeqComp(A\to A,\Gamma\Rightarrow \Delta)$, by the induction hypothesis, there exists a \textsc{contraction+cut}-free proof of $B, \Gamma_1,\Gamma_2\Rightarrow  \Delta'$. From this proof, we obtain a \textsc{contraction+cut}-free proof of $B\land C, \Gamma_1,\Gamma_2\Rightarrow  \Delta'$ as follows.
    \begin{prooftree}
    \def\fCenter{\Rightarrow}
    \Axiom$B, \Gamma_1,\Gamma_2\fCenter \Delta'$
    \RightLabel{{\color{red}{\quad $\land\sf{L}$}}}
    \UnaryInf$B\land C,\Gamma_1,\Gamma_2\fCenter\Delta'$
    \end{prooftree}
    Now note that, since $B\land C,\Gamma_1, \Gamma_2\Rightarrow \Delta'\preceq \Gamma\Rightarrow \Delta$, by \hyperlink{adm_neg_free}{\thmref{adm_neg_free}} it follows that $\Gamma\Rightarrow \Delta$ is obtainable from $B\land C, \Gamma_1,\Gamma_2\Rightarrow  \Delta'$ by suitable application(s) of \textsc{weakening} and/or \textsc{exchange}-rule(s).
    \end{subsubsubcase}

    \begin{subsubsubcase}
    $\Gamma'\Rightarrow \Delta'$ is the result of an application of ${\color{red}{\lor\sf{L}}}$. Then that part of the proof looks something like the following.
    \begin{prooftree}
    \def\fCenter{\Rightarrow}
    \AxiomC{$B,\Gamma^\ast\fCenter \Delta'$}
    \AxiomC{$C, \Gamma^\ast\fCenter \Delta'$}
    \RightLabel{{\color{red}{\quad $\lor\sf{L}$}}}
    \BinaryInf$B\lor C, \Gamma^\ast\fCenter\Delta'$
    \end{prooftree}Where $\Gamma'=B\lor C,\Gamma^\ast$.
    Now note that since $B\lor C, \Gamma^\ast\Rightarrow \Delta'$ has a \textsc{contraction+cut}-free proof, so does $B, \Gamma^\ast\Rightarrow \Delta'$ as well as $C, \Gamma^\ast\Rightarrow \Delta'$ (by \hyperlink{thm:invertibility}{\thmref{thm:invertibility}}). Since $A\to A$ was introduced earlier in the proof, it must be in the list $\Gamma^\ast$ (otherwise, refer to \hyperlink{subsubcase:2.2.2}{\subsubcaseref{subsubcase:2.2.2}}). Let $\Gamma^\ast:=\Gamma_1,A\to A, \Gamma_2$. Then $\Gamma'=B\lor C,\Gamma^\ast=B\lor C,\Gamma_1,A\to A, \Gamma_2$. Now, \begin{align*}&{\color{white}{\implies\,\,\,}}\Gamma'\Rightarrow \Delta'\preceq A\to A,\Gamma\Rightarrow \Delta\\&\implies B\lor C,\Gamma^\ast\Rightarrow \Delta'\preceq A\to A,\Gamma\Rightarrow \Delta\\&\implies B\lor C,\Gamma_1,A\to A, \Gamma_2\Rightarrow \Delta'\preceq A\to A,\Gamma\Rightarrow \Delta\\&\implies A\to A, B\lor C,\Gamma_1, \Gamma_2\Rightarrow \Delta'\preceq A\to A,\Gamma\Rightarrow \Delta\\&\implies B\lor C,\Gamma_1, \Gamma_2\Rightarrow \Delta'\preceq \Gamma\Rightarrow \Delta\end{align*}In other words,  $B\lor C,\Gamma_1,\Gamma_2\Rightarrow \Delta'$ is a subsequent of $\Gamma\Rightarrow \Delta$.
    
    Now, after the proofs of $B, \Gamma^\ast\Rightarrow \Delta'$ and $C, \Gamma^\ast\Rightarrow \Delta'$ we add the following lines, respectively.
    \begin{prooftree}
    \def\fCenter{\Rightarrow}
    \Axiom$B,\Gamma^\ast\fCenter \Delta'$
    \UnaryInf$B,\Gamma_1,A\to A, \Gamma_2\fCenter\Delta'$
    \RightLabel{{\color{red}{\quad $\sf{EL}$}}}
    \UnaryInf$A\to A,B,\Gamma_1,\Gamma_2\fCenter\Delta'$
    \end{prooftree}
    
    \begin{prooftree}
    \def\fCenter{\Rightarrow}
    \Axiom$C,\Gamma^\ast\fCenter \Delta'$
    \UnaryInf$C,\Gamma_1,A\to A, \Gamma_2\fCenter\Delta'$
    \RightLabel{{\color{red}{\quad $\sf{EL}$}}}
    \UnaryInf$A\to A,C,\Gamma_1,\Gamma_2\fCenter\Delta'$
    \end{prooftree}
    These constitute \textsc{contraction+cut}-free proofs of $A\to A,B,\Gamma_1,\Gamma_2\Rightarrow \Delta'$ and $A\to A,C,\Gamma_1,\Gamma_2\Rightarrow \Delta'$ respectively. Since $\SeqComp(A\to A, B, \Gamma_1,\Gamma_2\Rightarrow \Delta')<\SeqComp(A\to A,\Gamma\Rightarrow \Delta)$ and $\SeqComp(A\to A, C, \Gamma_1,\Gamma_2\Rightarrow \Delta')<\SeqComp(A\to A,\Gamma\Rightarrow \Delta)$, by the induction hypothesis we can conclude that there exist \textsc{contraction+cut}-free proofs of $B,\Gamma_1,\Gamma_2\Rightarrow  \Delta'$ and $C,\Gamma_1,\Gamma_2\Rightarrow  \Delta'$. From these proofs, we obtain a \textsc{contraction+cut}-free proof of $B\lor C,\Gamma_1,\Gamma_2\Rightarrow  \Delta'$ as follows.
    \begin{prooftree}
    \def\fCenter{\Rightarrow}
    \AxiomC{$B,\Gamma_1,\Gamma_2\fCenter\Delta'$}
    \AxiomC{$C,\Gamma_1,\Gamma_2\fCenter\Delta'$}
    \RightLabel{{\color{red}{\quad $\lor\sf{L}$}}}
    \BinaryInf$B\lor C,\Gamma_1,\Gamma_2\fCenter\Delta'$
    \end{prooftree}
    Now note that since $B\lor C,\Gamma_1,\Gamma_2\Rightarrow \Delta'$ is a subsequent of $\Gamma\Rightarrow \Delta$, by \hyperlink{adm_neg_free}{\thmref{adm_neg_free}} we can obtain a \textsc{contraction+cut}-free proof of $\Gamma\Rightarrow \Delta$ as well.
    \end{subsubsubcase}

    \begin{subsubsubcase}
    $\Gamma'\Rightarrow \Delta'$ is the result of an application of ${\color{red}{\neg\sf{L}}}$. Then that part of the proof looks something like the following.
    \begin{prooftree}
    \def\fCenter{\Rightarrow}
    \Axiom$\Gamma^\ast\fCenter \Delta', B$
    \RightLabel{{\color{red}{\quad $\neg\sf{L}$}}}
    \UnaryInf$\neg B, \Gamma^\ast\fCenter\Delta'$
    \end{prooftree}Where $\Gamma'=\neg B, \Gamma^\ast$.
    Now note that since $\neg B, \Gamma^\ast\Rightarrow \Delta'$ has a \textsc{contraction+cut}-free proof, so does $\Gamma^\ast\Rightarrow \Delta', B$ (by \hyperlink{thm:invertibility}{\thmref{thm:invertibility}}). Since $A\to A$ was introduced earlier in the proof, it must be in the list $\Gamma^\ast$ (otherwise, refer to \hyperlink{subsubcase:2.2.2}{\subsubcaseref{subsubcase:2.2.2}}). Without loss of generality, let $\Gamma^\ast:=\overline{\Gamma_1},A\to A,\overline{\Gamma_2}$. Then $\Gamma'=\neg B,\overline{\Gamma_1},A\to A,\overline{\Gamma_2}$. Now, \begin{align*}&{\color{white}{\implies\,\,\,}}\Gamma'\Rightarrow \Delta'\preceq A\to A,\Gamma\Rightarrow \Delta\\&\implies \neg B,\overline{\Gamma_1},A\to A,\overline{\Gamma_2}\Rightarrow \Delta'\preceq A\to A,\Gamma\Rightarrow \Delta\\&\implies A\to A,\neg B,\overline{\Gamma_1},\overline{\Gamma_2}\Rightarrow \Delta'\preceq A\to A,\Gamma\Rightarrow \Delta\\&\implies \neg B,\overline{\Gamma_1},\overline{\Gamma_2}\Rightarrow \Delta'\preceq \Gamma\Rightarrow \Delta\end{align*}In other words,  $\neg B,\overline{\Gamma_1},\overline{\Gamma_2}\Rightarrow \Delta'$ is a subsequent of $\Gamma\Rightarrow \Delta$.
    
    Now, we extend the proof of $\Gamma^\ast\Rightarrow \Delta', B$ as follows.
    \begin{prooftree}
    \def\fCenter{\Rightarrow}
    \Axiom$\Gamma^\ast\fCenter \Delta', B$
    \UnaryInf$\overline{\Gamma_1},A\to A,\overline{\Gamma_2}\fCenter\Delta', B$
    \RightLabel{{\color{red}{\quad $\sf{EL}$}}}
    \UnaryInf$A\to A,\overline{\Gamma_1},\overline{\Gamma_2}\fCenter\Delta', B$
    \end{prooftree}
    This constitutes a \textsc{contraction+cut}-free proof of $A\to A,\overline{\Gamma_1},\overline{\Gamma_2}\Rightarrow \Delta', B$. Since $\SeqComp(A\to A,\overline{\Gamma_1},\overline{\Gamma_2}\Rightarrow \Delta', B)< \SeqComp(A\to A,\Gamma\Rightarrow \Delta)$ by the induction hypothesis we can conclude that there exist a \textsc{contraction+cut}-free proof of $\overline{\Gamma_1},\overline{\Gamma_2}\Rightarrow \Delta', B$. From this proof, we obtain a \textsc{contraction+cut}-free proof of $\neg B,\overline{\Gamma_1},\overline{\Gamma_2}\Rightarrow \Delta'$ as follows.

    \begin{prooftree}
    \def\fCenter{\Rightarrow}
    \Axiom$\overline{\Gamma_1},\overline{\Gamma_2}\fCenter \Delta', B$
    \RightLabel{{\color{red}{\quad $\neg\sf{L}$}}}
    \UnaryInf$\neg B, \overline{\Gamma_1},\overline{\Gamma_2}\fCenter\Delta'$    \end{prooftree}
    
    Since $\neg B,\overline{\Gamma_1},\overline{\Gamma_2}\Rightarrow \Delta'$ is a subsequent of $\Gamma\Rightarrow \Delta$, applying \textsc{weakening} and \textsc{exchange} rules suitably, we can obtain $\Gamma\Rightarrow \Delta$, by \hyperlink{adm_neg_free}{\thmref{adm_neg_free}}.
    \end{subsubsubcase}

    \begin{subsubsubcase}
    $\Gamma'\Rightarrow \Delta'$ is the result of an application of ${\color{red}{\neg\sf{R}}}$. Then that part of the proof looks something like the following.
    \begin{prooftree}
    \def\fCenter{\Rightarrow}
    \Axiom$B, \Gamma'\fCenter \Delta^{\ast}$
    \RightLabel{{\color{red}{\quad $\neg\sf{R}$}}}
    \UnaryInf$\Gamma'\fCenter \Delta^\ast, \neg B$
    \end{prooftree}Where $\Delta'=\Delta^\ast, \neg B$.

    Now note that since $\Gamma'\Rightarrow \Delta^\ast, \neg B$ has a \textsc{contraction+cut}-free proof, so does $B, \Gamma'\Rightarrow \Delta^\ast$ (by \hyperlink{thm:invertibility}{\thmref{thm:invertibility}}). Since $A\to A$ was introduced earlier in the proof, it must be in the list $\Gamma'$ (otherwise, refer to \hyperlink{subsubcase:2.2.2}{\subsubcaseref{subsubcase:2.2.2}}). Let $\Gamma':=\Gamma_1,A\to A, \Gamma_2$. Now, \begin{align*}&{\color{white}{\implies\,\,\,}}\Gamma'\Rightarrow \Delta'\preceq A\to A,\Gamma\Rightarrow \Delta\\&\implies \Gamma'\Rightarrow \Delta^\ast, \neg B\preceq A\to A,\Gamma\Rightarrow \Delta\\&\implies \Gamma_1,A\to A, \Gamma_2\Rightarrow \Delta^\ast, \neg B\preceq A\to A,\Gamma\Rightarrow \Delta\\&\implies A\to A,\Gamma_1,\Gamma_2\Rightarrow \Delta^\ast, \neg B\preceq A\to A,\Gamma\Rightarrow \Delta\\&\implies \Gamma_1,\Gamma_2\Rightarrow \Delta^\ast, \neg B\preceq \Gamma\Rightarrow \Delta\end{align*}In other words,  $\Gamma_1,\Gamma_2\Rightarrow \Delta^\ast, \neg B$ is a subsequent of $\Gamma\Rightarrow \Delta$.

    Now, after the proof of $B, \Gamma'\Rightarrow \Delta^\ast$, we add the following steps.

    \begin{prooftree}
    \def\fCenter{\Rightarrow}
    \Axiom$B, \Gamma'\fCenter \Delta^\ast$
    \UnaryInf$B, \Gamma_1,A\to A, \Gamma_2\fCenter\Delta^\ast$
    \RightLabel{{\color{red}{\quad $\sf{EL}$}}}
    \UnaryInf$A\to A,B, \Gamma_1,\Gamma_2\fCenter\Delta^\ast$
    \end{prooftree}

    This constitutes a \textsc{contraction+cut}-free proof of $A\to A,B,\Gamma_1,\Gamma_2\Rightarrow \Delta^\ast$.  Since $\SeqComp(A\to A, B,\Gamma_1,\Gamma_2\Rightarrow \Delta^\ast)< \SeqComp(A\to A,\Gamma\Rightarrow \Delta)$, by the induction hypothesis, we can conclude that there exists a \textsc{contraction+cut}-free proof of $B,\Gamma_1,\Gamma_2\Rightarrow  \Delta^\ast$. From this proof, we obtain a \textsc{contraction+cut}-free proof of $\Gamma_1,\Gamma_2\Rightarrow  \Delta', \neg B$ as follows.

    \begin{prooftree}
    \def\fCenter{\Rightarrow}
    \Axiom$B,\Gamma_1,\Gamma_2\fCenter\Delta^\ast$
    \RightLabel{{\color{red}{\quad $\neg\sf{R}$}}}
    \UnaryInf$\Gamma_1,\Gamma_2\fCenter\Delta^\ast\neg B$
    \end{prooftree}
    Since $\Gamma_1,\Gamma_2\Rightarrow \Delta^\ast, \neg B$ is a subsequent of $\Gamma\Rightarrow \Delta$, by \hyperlink{adm_neg_free}{\thmref{adm_neg_free}} we can obtain $\Gamma\Rightarrow \Delta$.
    \end{subsubsubcase}

    \begin{subsubsubcase}
    $\Gamma\Rightarrow \Delta$ is the result of an application of ${\color{red}{\land\sf{R}}}$. Then that part of the proof looks something like the following.
    \begin{prooftree}
    \def\fCenter{\Rightarrow}
    \AxiomC{$\Gamma'\fCenter \Delta^\ast, B$}
    \AxiomC{$\Gamma'\fCenter \Delta^\ast, C$}
    \RightLabel{{\color{red}{\quad $\land\sf{R}$}}}
    \BinaryInfC{$\Gamma'\fCenter\Delta^\ast, B\land C$}
    \end{prooftree}Where $\Delta'=\Delta^\ast, B\land C$. Also, note that since $A\to A$ was introduced earlier in the proof, it must be in the list $\Gamma'$. Let $\Gamma':=\Gamma_1,A\to A,\Gamma_2$ (otherwise, refer to \hyperlink{subsubcase:2.2.2}{\subsubcaseref{subsubcase:2.2.2}}). Now, \begin{align*}&{\color{white}{\implies\,\,\,}}\Gamma'\Rightarrow \Delta'\preceq A\to A,\Gamma\Rightarrow \Delta\\&\implies \Gamma'\Rightarrow \Delta^\ast, B\land C\preceq A\to A,\Gamma\Rightarrow \Delta\\&\implies \Gamma_1,A\to A,\Gamma_2\Rightarrow \Delta^\ast, B\land C\preceq A\to A,\Gamma\Rightarrow \Delta\\&\implies A\to A,\Gamma_1,\Gamma_2\Rightarrow \Delta^\ast, B\land C\preceq A\to A,\Gamma\Rightarrow \Delta\\&\implies \Gamma_1,\Gamma_2\Rightarrow \Delta^\ast, B\land C\preceq \Gamma\Rightarrow \Delta\end{align*}In other words, $\Gamma_1,\Gamma_2\Rightarrow \Delta^\ast, B\land C$ is a subsequent of $\Gamma\Rightarrow \Delta$.

    Since $\Gamma'\Rightarrow \Delta^\ast, B\land C$ has a \textsc{contraction+cut}-free proof, so does $\Gamma'\Rightarrow \Delta^\ast, B$ and $\Gamma'\Rightarrow \Delta^\ast, C$ (by \hyperlink{thm:invertibility}{\thmref{thm:invertibility}}). We extend both these proofs respectively as follows.

    \begin{prooftree}
    \def\fCenter{\Rightarrow}
    \Axiom$\Gamma'\fCenter \Delta^\ast, B$
    \UnaryInf$\Gamma_1,A\to A,\Gamma_2\fCenter\Delta^\ast, B$
    \RightLabel{{\color{red}{\quad $\sf{EL}$}}}
    \UnaryInf$A\to A,\Gamma_1,\Gamma_2\fCenter\Delta^\ast, B$
    \end{prooftree}

    \begin{prooftree}
    \def\fCenter{\Rightarrow}
    \Axiom$\Gamma'\fCenter \Delta^\ast, C$
    \UnaryInf$\Gamma_1,A\to A,\Gamma_2\fCenter\Delta^\ast, C$
    \RightLabel{{\color{red}{\quad $\sf{EL}$}}}
   \UnaryInf$A\to A,\Gamma_1,\Gamma_2\fCenter\Delta^\ast, C$
    \end{prooftree}

    These now constitute \textsc{contraction+cut}-free proofs of $A\to A,\Gamma_1,\Gamma_2\Rightarrow \Delta^\ast, B$ and $A\to A,\Gamma_1,\Gamma_2\Rightarrow \Delta^\ast, C$ respectively. Since $\SeqComp(A\to A,\Gamma_1,\Gamma_2\Rightarrow \Delta^\ast, B)<\SeqComp(\neg A,\Gamma\Rightarrow \Delta)$ and $\SeqComp(A\to A,\Gamma_1,\Gamma_2\Rightarrow \Delta^\ast, C)< \SeqComp(A\to A,\Gamma\Rightarrow \Delta)$, by the induction hypothesis, there exist \textsc{contraction+cut}-free proofs of $ \Gamma_1,\Gamma_2\Rightarrow  \Delta', B$ and $ \Gamma_1,\Gamma_2\Rightarrow  \Delta', C$. From these proofs, we obtain a \textsc{contraction+cut}-free proof of $\Gamma_1,\Gamma_2\Rightarrow  \Delta^\ast, B\land C$ as follows.

    \begin{prooftree}
    \def\fCenter{\Rightarrow}
    \AxiomC{$\Gamma_1,\Gamma_2\fCenter \Delta^\ast, B$}
    \AxiomC{$\Gamma_1,\Gamma_2\fCenter \Delta^\ast, C$}
    \RightLabel{{\color{red}{\quad $\land\sf{R}$}}}
    \BinaryInfC{$\Gamma_1,\Gamma_2\fCenter\Delta^\ast, B\land C$}
    \end{prooftree}
    
    Since $\Gamma_1, \Gamma_2\Rightarrow \Delta^\ast, B\land C\preceq \Gamma\Rightarrow \Delta$ and since there is a \textsc{contraction+cut}-free proof of $\Gamma_1, \Gamma_2\Rightarrow \Delta^\ast, B\land C$, by \hyperlink{adm_neg_free}{\thmref{adm_neg_free}} it follows that $\Gamma\Rightarrow \Delta$ is obtainable from $\Gamma_1, \Gamma_2\Rightarrow \Delta^\ast, B\land C$ by suitable application(s) of \textsc{weakening} and/or \textsc{exchange}-rule(s).
    \end{subsubsubcase}

    \begin{subsubsubcase}
    $\Gamma'\Rightarrow \Delta'$ is the result of an application of ${\color{red}{\lor\sf{R}}}$. Then that part of the proof looks something like the following.
    \begin{prooftree}
    \def\fCenter{\Rightarrow}
    \Axiom$\Gamma'\fCenter \Delta^\ast, B$
    \RightLabel{{\color{red}{\quad $\lor\sf{R}$}}}
    \UnaryInf$\Gamma'\fCenter \Delta^\ast, B\lor C$
    \end{prooftree}Where $\Delta'=\Delta^\ast, B\lor C$.

    Since $\Gamma'\Rightarrow \Delta^\ast, B\lor C$ has a \textsc{contraction+cut}-free proof, so does $\Gamma'\Rightarrow \Delta^\ast, B$ (by \hyperlink{thm:invertibility}{\thmref{thm:invertibility}}). Since $A\to A$ was introduced earlier in the proof, it must be in the list $\Gamma'$ (otherwise, refer to \hyperlink{subsubcase:2.2.2}{\subsubcaseref{subsubcase:2.2.2}}). Let $\Gamma':=\Gamma_1,A\to A, \Gamma_2$. Now, \begin{align*}&{\color{white}{\implies\,\,\,}}\Gamma'\Rightarrow \Delta'\preceq A\to A,\Gamma\Rightarrow \Delta\\&\implies \Gamma'\Rightarrow \Delta^\ast, B\lor C\preceq A\to A,\Gamma\Rightarrow \Delta\\&\implies \Gamma_1,A\to A, \Gamma_2\Rightarrow \Delta^\ast, B\lor C\preceq A\to A,\Gamma\Rightarrow \Delta\\&\implies A\to A,\Gamma_1,\Gamma_2\Rightarrow \Delta^\ast, B\lor C\preceq A\to A,\Gamma\Rightarrow \Delta\\&\implies \Gamma_1,\Gamma_2\Rightarrow \Delta^\ast, B\lor C\preceq \Gamma\Rightarrow \Delta\end{align*}In other words,  $\Gamma_1,\Gamma_2\Rightarrow \Delta^\ast, B\lor C$ is a subsequent of $\Gamma\Rightarrow \Delta$.

    Now, after the proof of $\Gamma'\Rightarrow \Delta^\ast, B$, we add the following steps.

    \begin{prooftree}
    \def\fCenter{\Rightarrow}
    \Axiom$\Gamma'\fCenter \Delta^\ast, B$
    \UnaryInf$\Gamma_1,A\to A, \Gamma_2\fCenter\Delta^\ast, B$
    \RightLabel{{\color{red}{\quad $\sf{EL}$}}}
    \UnaryInf$A\to A,\Gamma_1,\Gamma_2\fCenter\Delta^\ast, B$
    \end{prooftree}

    This constitutes a \textsc{contraction+cut}-free proof of $A\to A,\Gamma_1,\Gamma_2\Rightarrow \Delta^\ast, B$.  Since we have $\SeqComp(A\to A,\Gamma_1,\Gamma_2\Rightarrow \Delta^\ast, B)< \SeqComp(A\to A,\Gamma\Rightarrow \Delta)$, by the induction hypothesis, we can conclude that there exists a \textsc{contraction+cut}-free proof of $\Gamma_1,\Gamma_2\Rightarrow  \Delta^\ast, B$. From this proof, we obtain a \textsc{contraction+cut}-free proof of $\Gamma_1,\Gamma_2\Rightarrow  \Delta^\ast, B\lor C$ as follows.

    \begin{prooftree}
    \def\fCenter{\Rightarrow}
    \Axiom$\Gamma_1,\Gamma_2\fCenter\Delta^\ast, B$
    \RightLabel{{\color{red}{\quad $\lor\sf{R}$}}}
    \UnaryInf$\Gamma_1,\Gamma_2\fCenter\Delta^\ast, B\lor C$
    \end{prooftree}
    Since $\Gamma_1,\Gamma_2\Rightarrow \Delta^\ast, B\lor C$ is a subsequent of $\Gamma\Rightarrow \Delta$, by \hyperlink{adm_neg_free}{\thmref{adm_neg_free}} we can obtain $\Gamma\Rightarrow \Delta$.
    \end{subsubsubcase}

    \begin{subsubsubcase}
    $\Gamma'\Rightarrow \Delta'$ is the result of an application of ${\color{red}{\to\sf{R}}}$. Then that part of the proof looks something like the following.
    \begin{prooftree}
    \def\fCenter{\Rightarrow}
    \Axiom$B,\Gamma'\fCenter \Delta^\ast, C$
    \RightLabel{{\color{red}{\quad $\to\sf{R}$}}}
    \UnaryInf$\Gamma'\fCenter\Delta^\ast, B\to C$
    \end{prooftree}i.e., in this case $\Delta'=\Delta^\ast, B\to C$.

    Since $\Gamma'\Rightarrow \Delta^\ast, B\to C$ has a \textsc{contraction+cut}-free proof, so does $B,\Gamma'\Rightarrow \Delta^\ast, C$ (by \hyperlink{thm:invertibility}{\thmref{thm:invertibility}}). Since $A\to A$ was introduced earlier in the proof, it must be in the list $\Gamma'$ (otherwise, refer to \hyperlink{subsubcase:2.2.2}{\subsubcaseref{subsubcase:2.2.2}}). Without loss of generality, let $\Gamma':=\Gamma_1,A\to A,\Gamma_2$. Now, \begin{align*}&{\color{white}{\implies\,\,\,}}\Gamma'\Rightarrow \Delta'\preceq A\to A,\Gamma\Rightarrow \Delta\\&\implies \Gamma_1,A\to A,\Gamma_2\Rightarrow \Delta^\ast, B\to C\preceq A\to A,\Gamma\Rightarrow \Delta\\&\implies A\to A,\Gamma_1,\Gamma_2\Rightarrow \Delta^\ast, B\to C\preceq A\to A,\Gamma\Rightarrow \Delta\\&\implies \Gamma_1,\Gamma_2\Rightarrow \Delta^\ast, B\to C\preceq \Gamma\Rightarrow \Delta\end{align*}In other words,  $\Gamma_1,\Gamma_2\Rightarrow \Delta^\ast, B\to C$ is a subsequent of $\Gamma\Rightarrow \Delta$.

    Now, after the proof of $B,\Gamma'\Rightarrow \Delta^\ast, C$, we add the following steps.
    \begin{prooftree}
    \def\fCenter{\Rightarrow}
    \Axiom$B,\Gamma'\fCenter\Delta^\ast, C$
    \UnaryInf$B,\Gamma_1, A\to A, \Gamma_2\fCenter\Delta^\ast, C$
    \RightLabel{{\color{red}{\quad $\sf{EL}$}}}
    \UnaryInf$A\to A, B,\Gamma_1,\Gamma_2\fCenter\Delta^\ast, C$
    \end{prooftree}

    This constitute a \textsc{contraction+cut}-free proof of $A\to A, B,\Gamma_1,\Gamma_2\Rightarrow \Delta^\ast, C$. Since we have $\SeqComp(A\to A,B,\Gamma_1,\Gamma_2\Rightarrow \Delta^\ast, C)<\SeqComp(A\to A,\Gamma\Rightarrow \Delta)$, by the induction hypothesis, we can conclude that there exists a \textsc{contraction+cut}-free proof of $B,\Gamma_1,\Gamma_2\Rightarrow \Delta^\ast, CA$. From this proof, we obtain a \textsc{contraction+cut}-free proof of $B,\Gamma_1,\Gamma_2\Rightarrow \Delta'$ as follows.

    \begin{prooftree}
    \def\fCenter{\Rightarrow}
    \Axiom$B,\Gamma_1,\Gamma_2\fCenter\Delta^\ast, C$
    \RightLabel{{\color{red}{\quad $\to\!\sf{R}$}}}
    \UnaryInf$\Gamma_1,\Gamma_2\fCenter\Delta^\ast, B\to C$
    \end{prooftree}
    Since $\Gamma_1,\Gamma_2\Rightarrow \Delta', B\to C$ is a subsequent of $\Gamma\Rightarrow \Delta$, applying \textsc{weakening} and \textsc{exchange} rules suitably, we can obtain $\Gamma\Rightarrow \Delta$, by \hyperlink{adm_neg_free}{\thmref{adm_neg_free}}.
    \end{subsubsubcase}    
\end{subsubcase}
\hypertarget{subsubcase:2.2.2}{\begin{subsubcase}{\label{subsubcase:2.2.2}}}
    If $A\to A$ does not occur in the sequent $\Gamma'\Rightarrow \Delta'$, we can immediately conclude that $\Gamma'\Rightarrow \Delta'$ is a subsequent of $\Gamma\Rightarrow \Delta$. Then, since $\Gamma'\Rightarrow \Delta'$ has a \textsc{contraction+cut}-free proof, by \hyperlink{adm_neg_free}{\thmref{adm_neg_free}},  $\Gamma\Rightarrow \Delta$ has a \textsc{contraction+cut}-free proof as well.
\end{subsubcase}
\end{subcase}
\end{case}
\end{proof}



\hypertarget{thm:con+cut_elim(IV)}{\begin{thm}[\textsc{Contraction+Cut-Elimination Theorem for \textbf{PLK} (Part III)}]{\label{thm:con+cut_elim(IV)}}Let $\Gamma\Rightarrow \Delta$ be a sequent. If $\Gamma\Rightarrow \Delta$ is provable in $\bf{PLK}$, then there is a \textsc{contraction+cut}-free $\bf{PLK}$-proof of $\Gamma\Rightarrow \Delta$. \end{thm}}

\begin{proof}
First, transform the given proof of $\Gamma\Rightarrow \Delta$ into a \textsc{contraction}-free proof by the procedure mentioned in \hyperlink{thm:con+cut_elim(I)}{\thmref{thm:con+cut_elim(I)}}. Choose any formula $A$ and obtain a \textsc{contraction}-free proof $A\to A, \Gamma\Rightarrow \Delta$ by an application of {\color{red}{$\mathsf{WL}$}}. Thus, we have obtained a \textsc{contraction}-free proof of $A\to A, \Gamma\Rightarrow \Delta$. Consequently, by \hyperlink{thm:con+cut_elim(II)}{\thmref{thm:con+cut_elim(II)}}, it follows that there exists a \textsc{contraction+cut}-free proof of $A\to A, \Gamma\Rightarrow \Delta$ as well. Finally, by \hyperlink{thm:con+cut_elim(III)}{\thmref{thm:con+cut_elim(III)}}, 
it follows that there exists a \textsc{contraction+cut}-free proof of $\Gamma\Rightarrow \Delta$. 
\end{proof}

\hypertarget{rem:additive_choice}{\begin{rem}{\label{rem:additive_choice}}\hyperlink{thm:con+cut_elim(I)}{\thmref{thm:con+cut_elim(I)}} has been one of the crucial ingredients for proving the above theorem. It should be pointed out that if one considered the multiplicative version of \textsc{cut} instead of the additive version, the result would not hold (see \cite{Kashima1997} for details). Furthermore, the $\SeqComp$ function, so crucial in proving \hyperlink{thm:con+cut_elim(III)}{\thmref{thm:con+cut_elim(III)}}, induces an ordering in the rules; namely, the $\SeqComp$-value of every premise of a rule is less than or equal to the value of the conclusion. This order is essentially reversed in the case of (additive) \textsc{cut} and \textsc{contraction}-rules. To sum up, in each case, we could conclude something about the ordering of the premises and conclusion in a \textit{uniform} manner. This is not the case if we consider $M\textsc{cut}$. For example, if $\Gamma,\Delta, \Pi,\Phi$ are all nonempty, having at least one non-atomic formula, and if $p$ is a propositional variable, then considering the following instance of $M\textsc{cut}$,
 \begin{prooftree}
    \def\fCenter{\Rightarrow}
    \Axiom$\Gamma\fCenter\Delta, \neg p$
    \Axiom$\neg p, \Pi\fCenter\Phi$
    \RightLabel{{\color{red}{\quad $\textsc{cut}$}}}
    \BinaryInf$\Gamma,\Pi\fCenter\Delta,\Phi$
\end{prooftree}
we note that the $\SeqComp$-value of each premise is strictly less than that of the conclusion. If, on the other hand, we consider the following instance of $M\textsc{cut}$ where each element of $\Gamma,\Delta, \Pi, \Phi$ is a propositional variable,
 \begin{prooftree}
    \def\fCenter{\Rightarrow}
    \Axiom$\Gamma\fCenter\Delta, p$
    \Axiom$p, \Pi\fCenter\Phi$
    \RightLabel{{\color{red}{\quad $\textsc{cut}$}}}
    \BinaryInf$\Gamma,\Pi\fCenter\Delta,\Phi$
\end{prooftree} we note that the $\SeqComp$-value of each premise is the same as that of the conclusion. If, on the other hand, we consider the following instance of $M\textsc{cut}$: 
 \begin{prooftree}
    \def\fCenter{\Rightarrow}
    \Axiom$\Gamma\fCenter\Delta, A$
    \Axiom$A, \Pi\fCenter\Phi$
    \RightLabel{{\color{red}{\quad $\textsc{cut}$}}}
    \BinaryInf$\Gamma,\Pi\fCenter\Delta,\Phi$
\end{prooftree} where $\mathsf{Comp}(A)>\SeqComp(\Gamma,\Pi\Rightarrow \Delta,\Phi)$, it follows that the value of $\SeqComp$ for each of the premise is strictly greater than that of the conclusion. Consequently, no uniform conclusion about the ordering of rules via the number of connectives can be maintained for the $M\textsc{cut}$. Although this doesn't exclude the possibility of finding some \textit{other} function relative to which the rules can be ordered, it is certain that \textit{that} function cannot be $\SeqComp$. 
These are the main reasons for choosing the additive version of \textsc{cut} instead of the multiplicative one.
\end{rem}}

\subsection{Non-Eliminability of the Other Rules of \textbf{PLK}}
So far, we have been successful in eliminating the \textsc{contraction}-rules and the \textsc{cut}-rule. Is
there any other rule in \textbf{PLK} which can be eliminated? In this section, we delve into this question and show that no other \textbf{PLK}-rule can be eliminated. In fact, we show something stronger. We show that for each rule except \textsc{contraction}-rules and the \textsc{cut}-rule there exists a sequent, the proof of which can't be done without using that rule. It is important to note that all of these conclusions follow as corollaries of \hyperlink{thm:con+cut_elim(IV)}{\thmref{thm:con+cut_elim(IV)}}. 

\hypertarget{thm:non_elim_struct}{\begin{thm}[\textsc{Non-Eliminability of the Structural Rules of \textbf{PLK}}]{\label{thm:non_elim_struct}} Let $p$ and $q$ denote distinct propositional variables. Then, the following statements hold.
\begin{enumerate}[label=(\arabic*)]
    \item $q,p\Rightarrow p $ is not provable without using $\color{red}{\mathsf{WL}}$.
    \item $p\Rightarrow p,q$ is not provable without using $\color{red}{\mathsf{WR}}$. 
    \item $p\Rightarrow q,p$ is not provable without using $\color{red}{\mathsf{EL}}$.
    \item $p,q\Rightarrow p$ is not provable without using $\color{red}{\mathsf{ER}}$.
\end{enumerate}\end{thm}}

\begin{proof}
We prove (1) only. The proofs for the other rules are similar.

Note that $q,p\Rightarrow p $ is a \textbf{PLK}-provable sequent as shown below: 
\begin{prooftree}
    \def\fCenter{\Rightarrow}
    \Axiom$p\fCenter p$
    \RightLabel{{\color{red}{\quad $\sf{WL}$}}}
    \UnaryInf$q,p\fCenter p$
\end{prooftree}
Consider any \textbf{PLK}-proof of $q,p\Rightarrow p$. Then, by \hyperlink{thm:con+cut_elim(IV)}{\thmref{thm:con+cut_elim(IV)}}, there exists a \textsc{contraction+cut}-free proof of the same. Consider any such proof. Note that $q,p\Rightarrow p$ can't be the conclusion of any instance of any logical rule since it does not contain any connective. So, it must be the conclusion of an instance of some \textsc{exchange}-rule or \textsc{weakening}-rule.

But, it can't be the conclusion of any instance of $\color{red}{\mathsf{WR}}$ because that would imply the provability of the sequent $q,p\Rightarrow \Lambda$, an impossibility due to \hyperlink{thm:PLK_cons}{\thmref{thm:PLK_cons}}. Also, it can't be the conclusion of any instance of $\color{red}{\mathsf{ER}}$ because to apply $\color{red}{\mathsf{ER}}$, one needs at least two wffs on the right of the premise.

Hence, $q,p\Rightarrow p$ can only be the conclusion of some instance of $\color{red}{\mathsf{WL}}$ or $\color{red}{\mathsf{EL}}$. In case it is the conclusion of an instance of $\color{red}{\mathsf{WL}}$, it must be the following:

\begin{prooftree}
    \def\fCenter{\Rightarrow}
    \Axiom$p\fCenter p$
    \RightLabel{{\color{red}{\quad $\sf{WL}$}}}
    \UnaryInf$q,p\fCenter p$
\end{prooftree}In case it is not the conclusion of an instance of $\color{red}{\mathsf{WL}}$, it must be the conclusion of an instance of  $\color{red}{\mathsf{EL}}$. The last two lines of the proof then look like the following:

\begin{prooftree}
    \def\fCenter{\Rightarrow}
    \Axiom$p,q\fCenter p$
    \RightLabel{{\color{red}{\quad $\sf{EL}$}}}
    \UnaryInf$q,p\fCenter p$
\end{prooftree}
Consider now the sequent $p,q\Rightarrow p$. By an argument similar to the above, we can conclude that it can't be the conclusion of any instance of \textsc{cut}, \textsc{contraction}, $\color{red}{\mathsf{WR}}$, $\color{red}{\mathsf{ER}}$ or logical rules. However, in this case, it can't be the conclusion of an instance $\color{red}{\mathsf{WL}}$ either since that would imply the \textbf{PLK}-provability of $q\Rightarrow p$, an impossibility $-$ by \hyperlink{thm:PLK_cons}{\thmref{thm:PLK_cons}} $-$ since $p\ne q$. Hence, it must be the conclusion of an instance of $\color{red}{\mathsf{EL}}$. Consequently, the last two lines of the proof will look something like the following:

\begin{prooftree}
    \def\fCenter{\Rightarrow}
    \Axiom$q,p\fCenter p$
    \RightLabel{{\color{red}{\quad $\sf{EL}$}}}
    \UnaryInf$p,q\fCenter p$
\end{prooftree}
and we are back to the same sequent we started with. So, we continue the argument and go upwards in the proof, i.e., we check whether the concerned sequent is the conclusion of some instance of an \textsc{exchange}-rule. At the first occurrence of a sequent, which is not the conclusion of any instance of an \textsc{exchange}-rule, we stop (such a sequent must exist since the proof is \textit{finite}). Hence, this sequent must be the conclusion of an instance of $\color{red}{\mathsf{WL}}$, and we are done.   
\end{proof}

\hypertarget{thm:non_elim_logic}{\begin{thm}[\textsc{Non-Eliminability of the Logical Rules of \textbf{PLK}}]{\label{thm:non_elim_logic}}Let $p$ and $q$ denote distinct propositional variables. Then, the following statements hold.
\begin{enumerate}[label=(\arabic*)]
    \item $\neg p,p\Rightarrow \Lambda $ is not provable without using $\color{red}{\neg\mathsf{L}}$.
    \item $\Lambda\Rightarrow p,\neg p$ is not provable without using $\color{red}{\neg\mathsf{R}}$.
    \item $p\land q\Rightarrow p$ is not provable without using $\color{red}{\land\mathsf{L}}$. 
    \item $p\Rightarrow p\land p$ is not provable without using $\color{red}{\land\mathsf{R}}$.
    \item $p\lor p\Rightarrow p$ is not provable without using $\color{red}{\lor\mathsf{L}}$. 
    \item $p\Rightarrow p\lor q$ is not provable without using $\color{red}{\lor\mathsf{R}}$.
    \item $p\to q,p\Rightarrow q$ is not provable without using $\color{red}{\to\mathsf{L}}$.
    \item $\Lambda\Rightarrow p\to p$ is not provable without using $\color{red}{\to\mathsf{R}}$.
\end{enumerate}\end{thm}}

\begin{proof}
We proof only (7). The proofs for other rules are similar. 

Note that $p\to q,p\Rightarrow q$ is a \textbf{PLK}-provable sequent as shown below. 
\begin{prooftree}
    \def\fCenter{\Rightarrow}
    \Axiom$p\fCenter p$
    \LeftLabel{{\color{red}{$\sf{WR}$\quad }}}
    \UnaryInf$p\fCenter p,q$
    \LeftLabel{{\color{red}{$\sf{ER}$\quad }}}
    \UnaryInf$p\fCenter q,p$
    \Axiom$q\fCenter q$
    \RightLabel{{\color{red}{\quad $\sf{WL}$}}}
    \UnaryInf$p,q\fCenter q$
    \RightLabel{{\color{red}{\quad $\sf{EL}$}}}
    \UnaryInf$q,p\fCenter q$
    \RightLabel{{\color{red}{\quad $\to\sf{L}$}}}
    \BinaryInf$p\to q,p\fCenter q$
\end{prooftree}
Consider any \textbf{PLK}-proof of $p\to q,p\Rightarrow q $. Then, by \hyperlink{thm:con+cut_elim(IV)}{\thmref{thm:con+cut_elim(IV)}}, there exists a \textsc{contraction+cut}-free proof of the same. Consider any such proof. Note that $p\to q,p\Rightarrow q $ can't be the conclusion of any instance of any logical rule other than ${\color{red}{\to\sf{L}}}$ since it does not contain any other connective. It can't be the conclusion of any instance of ${\color{red}{\sf{WL}}}$ or ${\color{red}{\sf{WR}}}$ because in the former case we would be able to conclude that $p\Rightarrow q$ is \textbf{PLK}-provable, while in the latter we will be able to conclude that $p\to q, p\Rightarrow \emptyset$ is \textbf{PLK}-provable. The first case is impossible by \hyperlink{thm:PLK_cons}{\thmref{thm:PLK_cons}} (since $p$ and $q$ are distinct variables); the second case is also impossible since provability of $p\to q, p\Rightarrow \emptyset$ will imply (due to \hyperlink{thm:invertibility}{\thmref{thm:invertibility}}) that there is a proof of $q,p\Rightarrow  \emptyset$, contrary to \hyperlink{thm:PLK_cons}{\thmref{thm:PLK_cons}}.

Furthermore, $p\to q,p\Rightarrow q $ can't be the conclusion of any instance of $\color{red}{\mathsf{ER}}$ because to apply $\color{red}{\mathsf{ER}}$ one needs at least two wffs on the right in the premise.

Hence, $p\to q,p\Rightarrow q$ can only be the conclusion of some instance of $\color{red}{\to\mathsf{L}}$ or $\color{red}{\mathsf{EL}}$. Without loss of generality, we may assume that $p\to q,p\Rightarrow q$ is the conclusion of an instance of $\color{red}{\to\mathsf{L}}$ because if it is the conclusion of an instance of $\color{red}{\mathsf{EL}}$, we will go upwards in the proof to find the sequent which is not the conclusion of an instance of $\color{red}{\mathsf{EL}}$. This sequent must be the conclusion of an instance of $\color{red}{\to\mathsf{L}}$, and in either case, that part of the proof must look something like the following:

\begin{prooftree}
    \def\fCenter{\Rightarrow}
    \Axiom$p\fCenter q,p$
    \Axiom$q,p\fCenter q$
    \RightLabel{{\color{red}{\quad $\to\sf{L}$}}}
    \BinaryInf$p\to q,p\fCenter q$
\end{prooftree}
Consider now the sequent $q,p\Rightarrow q$. By an argument similar to the above, we will be able to conclude that it can't be the conclusion of any instance of \textsc{cut}, \textsc{contraction}, $\color{red}{\mathsf{WR}}$, $\color{red}{\mathsf{ER}}$ or logical rules. However, in this case, it can't be the conclusion of any instance of $\color{red}{\mathsf{WL}}$ either since that would imply the \textbf{PLK}-provability of $p\Rightarrow q$, an impossibility (again by \hyperlink{thm:PLK_cons}{\thmref{thm:PLK_cons}}). Hence, it must be the conclusion of some instance of $\color{red}{\mathsf{EL}}$. Consequently, the last two lines of the proof of $q,p\Rightarrow q$ in the proof of $p\to q,p\Rightarrow q$ will look something like the following:

\begin{prooftree}
    \def\fCenter{\Rightarrow}
    \Axiom$p,q\fCenter q$
    \RightLabel{{\color{red}{\quad $\sf{EL}$}}}
    \UnaryInf$q,p\fCenter q$
\end{prooftree}
We continue the argument and go upwards in the proof, i.e., we check whether the corresponding sequent is the conclusion of some instance of \textsc{exchange}-rules. At the first occurrence of a sequent, which is not the conclusion of any instance of \textsc{exchange}-rules, we stop (such a sequent must exist since the proof is \textit{finite}). Hence, this sequent must be the conclusion of an instance of $\color{red}{\mathsf{WL}}$, and it must be $p,q\Rightarrow q$. By an analogous argument we can show that the case for the sequent $p\Rightarrow q, p$ will be similar. To sum up, the proof of $p\to q,p\Rightarrow q$ will always look something like the following: 
\begin{prooftree}
    \def\fCenter{\Rightarrow}
    \Axiom$p\fCenter p$
    \LeftLabel{{\color{red}{$\sf{WR}$\quad }}}
    \UnaryInf$p\fCenter p,q$
    \LeftLabel{{\color{red}{$\sf{ER}$\quad }}}
     \UnaryInf$p\fCenter q,p$
    \LeftLabel{{\color{red}{$\sf{ER}$\quad }}}
    \def\fCenter{\vdots}
    \UnaryInf$\phantom{p}\fCenter\phantom{p}$
    \LeftLabel{{\color{red}{$\sf{ER}$\quad }}}
    \def\fCenter{\Rightarrow}
    \UnaryInf$p\fCenter q,p$
    \Axiom$q\fCenter q$
    \RightLabel{{\color{red}{\quad $\sf{WL}$}}}
    \UnaryInf$p,q\fCenter q$
    \RightLabel{{\color{red}{\quad $\sf{EL}$}}}
     \def\fCenter{\vdots}
    \UnaryInf$\phantom{p}\fCenter\phantom{p}$
    \RightLabel{{\color{red}{\quad $\sf{EL}$}}}
    \def\fCenter{\Rightarrow}
    \UnaryInf$q,p\fCenter q$
    \RightLabel{{\color{red}{\quad $\to\sf{L}$}}}
    \BinaryInf$p\to q,p\fCenter q$
    \def\fCenter{\vdots}
    \UnaryInf$\phantom{p}\fCenter\phantom{p}$
    \def\fCenter{\Rightarrow}
    \UnaryInf$p\to q,p\fCenter q$
\end{prooftree}
and we are done.   

\end{proof}

\subsection{Observations}
We now make a couple of observations.
\begin{enumerate}[label=(\alph*)]
\item The proofs of the \textsc{cut}-elimination theorem, detailed in the previous sections, cannot be applied in their present form to prove \textsc{cut}-elimination theorem for any system in which either one or both of the \textsc{weakening}-rules and/or the \textsc{exchange}-rules are restricted. This, however, does not restrict the possibility of eliminating some restricted form of \textsc{cut} from the sequent system by applying the same technique. 

\item As remarked earlier, in each case of the proof of \hyperlink{ce1}{\thmref{ce1}}, the only thing which makes the induction work, apart from \hyperlink{ml1}{\thmref{ml1}}, is the fact that the complexity of each of the premise sequent of a logical rule is strictly less than that of the conclusion sequent. In particular, it should be noted that in the argument, the number of premises is also not important. 

\item We also note that since the notion of sequent complexity is independent of particular connectives, the proof can be easily adapted to other sequent systems as well. However, analysing carefully, we find that even the notion of complexity is not necessary, and we can generalise the argument even more. 

\item The only structural rules needed were \textsc{weakening} and \textsc{exchange}.

\item The fact that the sequent $\Lambda\Rightarrow \Lambda$ has no proof in \textbf{PLK} (proven in \hyperlink{thm:PLK_cons}{\thmref{thm:PLK_cons}}) is crucial for the proof of  \hyperlink{thm:con+cut_elim(III)}{\thmref{thm:con+cut_elim(III)}}.
\end{enumerate}

These observations motivate us to generalise the argument to a broader class of sequent systems. In the next section, we will give a sufficient condition for a sequent system to have the property of \textsc{cut}-elimination.

\section{Normal Sequent Structures}

In this section, we generalise the arguments of the previous section to a broader class of sequent systems, taking note of the observations just mentioned. We introduce the notion of \textit{normal sequent structures} and prove two \textsc{rule}-elimination theorems for the same. 

We begin with the following definitions.

\begin{defn}[\textsc{Sequent}]Let $\mathscr{L}$ be a set and let $\Rightarrow\subseteq \boldsymbol{FO}(\mathscr{L})\times\boldsymbol{FO}(\mathscr{L})$, where $\boldsymbol{FO}(\mathscr{L})$ denotes the set of all finite ordered lists of elements of $\mathscr{L}$. Then $(\Gamma,\Delta)\in \Rightarrow$ will be said to be a \textit{$\mathscr{L}$-sequent with respect to $\Rightarrow$}. In this case, instead of writing $(\Gamma,\Delta)\in \Rightarrow$, we will write $\Gamma\Rightarrow \Delta$. The \textit{set of all $\mathscr{L}$-sequents with respect to $\Rightarrow$}, will henceforth be denoted by $\mathsf{Seq}(\mathscr{L},\Rightarrow )$.%
\end{defn}

\hypertarget{def:rule}{\begin{defn}[\textsc{Rule}]{\label{def:rule}}Let $\mathsf{Seq}(\mathscr{L},\Rightarrow )$ be the set of all $\mathscr{L}$-sequents with respect to $\Rightarrow$. A set $R$ is said to be a \textit{$\mathsf{Seq}(\mathscr{L},\Rightarrow )$-rule} if, $$\emptyset \subsetneq R\subseteq \mathsf{Seq}(\mathscr{L},\Rightarrow )^n\times \mathsf{Seq}(\mathscr{L},\Rightarrow )$$for some $n\in \mathbb{N}~(\text{where}~\mathbb{N}:=\{0,1,\ldots\})$. It is said to be a \textit{$\mathsf{Seq}(\mathscr{L},\Rightarrow )$-axiomatic rule} if $\emptyset \subsetneq R\subseteq \{\emptyset\}\times \mathsf{Seq}(\mathscr{L},\Rightarrow )$, else, a \textit{non-axiomatic rule}. An \textit{instance} of a rule $R$ is an element of $R$. Given $((\Gamma_1\Rightarrow \Delta_1,\ldots, \Gamma_n\Rightarrow \Delta_n),\Gamma\Rightarrow \Delta)\in R$, $\Gamma_i\Rightarrow \Delta_i$ is to be a \textit{premise} of this instance of $R$ for each $i\in \{1,\ldots,n\}$ and $\Gamma\Rightarrow \Delta$ is said to be a \textit{conclusion} of this instance of $R$. In this case, we will also sometimes say $\Gamma\Rightarrow \Delta$ to be \textit{obtained by an application of $R$}.\end{defn}}

Recall that in \hyperlink{sec:cut_elim(PLK)}{\secref{sec:cut_elim(PLK)}}, the language for \textbf{PLK} is defined to be an algebra, constructed recursively over a countably infinite set of propositional variables using a finite set of connectives. In this section, we generalise this and consider our language to be a formula algebra. To be more precise, we assume that $\mathscr{L}$ is the set of \textit{formulas} defined inductively in the usual way over a nonempty set of variables $V$ using a finite set of connectives/operators called the \textit{signature/type}. The formula algebra has the universal mapping property for the class of all algebras of the same type as $\mathscr{L}$ over $V$, i.e., any function $f : \mathscr{P}\to A$, where $A$ is the universe of an algebra \textbf{A} of
the same type as $\mathscr{L}$, can be uniquely extended to a homomorphism from $\mathscr{L}$ to $\bf{A}$ (see \cite{Font2016} for more details).

\begin{defn}[\textsc{Normal Sequent Structure}]{\label{pnseqsys}}
Let $\mathscr{L}$ be a formula algebra of some type over a set of  variables $V$ and $\mathsf{Seq}(\mathscr{L},\Rightarrow )$ be the set of all $\mathscr{L}$-sequents with respect to $\Rightarrow$ and $\vdash\,\subseteq \mathcal{P}(\mathsf{Seq}(\mathscr{L},\Rightarrow ))\times \mathsf{Seq}(\mathscr{L},\Rightarrow )$. The pair $(\mathsf{Seq}(\mathscr{L},\Rightarrow ),\vdash)$ is then said to be a \textit{normal sequent structure} if the following conditions hold.
\begin{enumerate}[label=(\roman*)]
\item $p\Rightarrow  p$ for all $p\in V$. Furthermore, $\emptyset\vdash p\Rightarrow  p$ for all $p\in V$.
\item Both the \textsc{weakening} and \textsc{exchange} rules are included in the set of rules. 
\hypertarget{pnseqsys(3)}{\item{\label{pnseqsys(3)}} $\emptyset\vdash \Gamma\Rightarrow  \Delta$ means that there is a proof of $\Gamma\Rightarrow \Delta$; or, alternatively, that $\Gamma\Rightarrow  \Delta$ is \textit{provable}. Proofs are, as in \hyperlink{sec:cut_elim(PLK)}{\secref{sec:cut_elim(PLK)}}, built as trees of sequents (i.e. with nodes labelled with sequents), with leaves being conclusions of an instance of some axiomatic rule and the proven sequent as a root. Formally this notion may be defined inductively as follows:
\begin{enumerate}[label=$\bullet$]
\item If $S_0$ is the conclusion of an instance of an axiomatic rule, then it is a proof of a sequent $S_0$.
\item If $\mathcal{D}_1$ is a proof of a sequent $S_0$, then $\dfrac{\mathcal{D}_1}{S_1}$ is a proof of a sequent $S_1$
, provided that $S_0$ is an instance of the premise and $S_1$ an instance of the conclusion of some one-premise rule.
\item If $\mathcal{D}_1$ is a proof of a sequent $S_0$ and $\mathcal{D}_2$ is a proof of a sequent $S_1$, then
$\dfrac{\mathcal{D}_1~~~\mathcal{D}_2}{S_2}$ is a proof of a sequent $S_2$, provided that $S_0$ and $S_1$ are instances of
the premises and $S_2$ an instance of the conclusion of some two-premise rule.
\item If $\mathcal{D}_1$ is a proof of a sequent $S_0$; $\mathcal{D}_2$ is a proof of a sequent $S_1$ and $\mathcal{D}_3$ is a proof of a sequent $S_2$, then
$\dfrac{\mathcal{D}_1~~~\mathcal{D}_2~~~\mathcal{D}_3}{S_3}$ is a proof of a sequent $S_33$, provided that $S_0,S_1$ and $S_2$ are instances of the premises, and $S_3$ an instance of the conclusion of some three-premise rule.
\end{enumerate}
... and so on, for each $n\in \mathbb{N}\setminus\{0\}$ and for each $n$-premise rule.}
\hypertarget{pnseqsys(4)}{\item{\label{pnseqsys(4)}} If $\Gamma\Rightarrow \Delta$ is provable in $\mathsf{Seq}(\mathscr{L},\Rightarrow )$ and if $\Gamma$ and $\Delta$ both contains only variables, then there exists $p\in V$ that is common to both the lists $\Gamma$ and $\Delta$.}
\hypertarget{pnseqsys(5)}{\item{\label{pnseqsys(5)}} If $\Gamma\Rightarrow \Delta$ is provable in $\mathsf{Seq}(\mathscr{L},\Rightarrow )$, then either $\Gamma$ is not the empty list or $\Delta$ is not the empty list.}
\end{enumerate}
\end{defn}

In the case of \textbf{PLK}, $\SeqComp$ satisfies the following interesting property: for each rule, if the side formulae are non-atomic, then (except \textsc{cut} and \textsc{contraction}-rules) the $\SeqComp$-value of each of the premise is less than or equal to that of the conclusion. This is not true for \textsc{cut} and \textsc{contraction}-rules, and they were eliminable. The motivation for the following definition was the (conjectural) belief that this is not a coincidence. To be more precise, we think that both the \textsc{cut} and \textsc{contraction}-rules were eliminable \textit{precisely because} they did not harmonise with other rules. $\SeqComp$ merely `witnessed' this. 

\hypertarget{witfun}{\begin{defn}[\textsc{Witness Function}]{\label{witfun}}
Let $(\mathsf{Seq}(\mathscr{L},\Rightarrow ),\vdash)$ be a normal sequent structure and $v:\mathsf{Seq}(\mathscr{L},\Rightarrow )\to\mathbb{N}$. Then $v$ is said to be a \textit{witness function for $(\mathsf{Seq}(\mathscr{L},\Rightarrow ),\vdash)$ relative to the rule $T$} if the following conditions hold.
\begin{enumerate}[label=(\roman*)]
    \hypertarget{witfun(1)}{\item{\label{witfun(1)}} There exists $n_0\in \mathbb{N}$ such that if for some sequent $\Gamma\Rightarrow  \Delta$, $v(\Gamma\Rightarrow  \Delta)=n_0$ then both $\Gamma$ and $\Delta$ contain at most variables. Furthermore, for all $\Gamma\Rightarrow \Delta\in \mathsf{Seq}(\mathscr{L},\Rightarrow )$, $v(\Gamma\Rightarrow \Delta)\ge n_0$.}
    \hypertarget{witfun(2)}{\item{\label{witfun(2)}} Suppose $R$ is a $\mathsf{Seq}(\mathscr{L},\Rightarrow )$-non-axiomatic rule. We say that $v$ \emph{respects} (respectively, \emph{strictly respects}) an instance $((\Gamma_1\Rightarrow \Delta_1,\ldots,\Gamma_m\Rightarrow \Delta_m), \Gamma\Rightarrow \Delta)$ of $R$ if $v(\Gamma_i\Rightarrow \Delta_i)\le v(\Gamma\Rightarrow \Delta)$ (respectively, if $v(\Gamma_i\Rightarrow \Delta_i)< v(\Gamma\Rightarrow \Delta)$), for all $i\in \{1,\ldots,m\}$.}
    \item $v$ does not respect any instance of $T$. 
\end{enumerate}
\end{defn}}

\hypertarget{nseqsyscutelim1}{\begin{thm}[\textsc{Rule-Elimination Theorem for NSSs (I)}]{\label{nseqsyscutelim1}}
Let $(\mathsf{Seq}(\mathscr{L},\Rightarrow ),\vdash)$ be a normal sequent structure having a witness function relative to a rule $T$ (which is disjoint from \textsc{weakening} and \textsc{exchange} rules). Suppose that for every $\mathsf{Seq}(\mathscr{L},\Rightarrow )$-provable sequent $\Gamma\Rightarrow\Delta$, the following statements hold.
\begin{enumerate}[label=(\arabic*)]
    \item If $n_0<v(\Gamma\Rightarrow \Delta)$ then there exists a provable $\Gamma'\Rightarrow \Delta'$ such that $v(\Gamma'\Rightarrow \Delta')<v(\Gamma\Rightarrow \Delta)$.
    \item If $\Gamma'\Rightarrow \Delta'$ has a $T$-free proof, then so does $\Gamma\Rightarrow \Delta$.
\end{enumerate}
Then every provable sequent has a $T$-free proof. 
\end{thm}}

\begin{proof}
Let $v:\mathsf{Seq}(\mathscr{L},\Rightarrow )\to \mathbb{N}$ be the witness function relative to $T$ and $\Gamma\Rightarrow \Delta$ be a provable sequent. We will show this by induction on the values of $v$. 

$\\$\underline{\textsc{Base Case:}} In the base case, $n = n_0$; thus every formula in the sequent $\Gamma\Rightarrow \Delta$ is a propositional variable (by \hyperlink{witfun(1)}{\defref{witfun}(i)}). Since the sequent is provable, there must be some variable $p$, which occurs both in $\Gamma$ and $\Delta$ (by \hyperlink{pnseqsys(4)}{\defref{pnseqsys}(iv))}). Since $p\Rightarrow p$ is provable, $\Gamma\Rightarrow \Delta$ can be proved without using $T$ by application (or applications, as needed)  of the \textsc{weakening} and/or \textsc{exchange} rules.%
\footnote{Here we use the fact that $T$ is disjoint from both \textsc{weakening} and \textsc{exchange}.}

$\\$\underline{\textsc{Induction Hypothesis:}} For every sequent $\Gamma\Rightarrow \Delta$ such that $n_0\le v(\Gamma\Rightarrow \Delta)<n$, if there is a $\mathsf{Seq}(\mathscr{L},\Rightarrow )$-proof of $\Gamma\Rightarrow \Delta$ , then there is a $T$-free $\mathsf{Seq}(\mathscr{L},\Rightarrow )$-proof of $\Gamma\Rightarrow \Delta$.

$\\$\underline{\textsc{Induction Step:}} Let $\Gamma\Rightarrow \Delta$ be a sequent such that $v(\Gamma\Rightarrow \Delta)=n$. Suppose also that $\Gamma\Rightarrow \Delta$ is $\mathsf{Seq}(\mathscr{L},\Rightarrow )$-provable. Then, by (1), it follows that there exists a provable $\Gamma'\Rightarrow \Delta'$ with $v(\Gamma'\Rightarrow \Delta')<v(\Gamma\Rightarrow \Delta)$ such that $\Gamma'\Rightarrow \Delta'$ is provable. Since $v(\Gamma'\Rightarrow \Delta')<v(\Gamma\Rightarrow \Delta)$ by the induction hypothesis it follows that $\Gamma'\Rightarrow \Delta'$ has a $T$-free proof. Consequently, by (2), it follows that $\Gamma\Rightarrow \Delta$ has a $T$-free proof as well.
\end{proof}

\begin{rem}
Note that in this proof the fact that $\mathscr{L}$ is an algebra is needed only to prove the base case.
\end{rem}

\hypertarget{nseqsyscutelim2}{\begin{thm}[\textsc{Rule-Elimination Theorem for NSSs (II)}]{\label{nseqsyscutelim2}}
Let $(\mathsf{Seq}(\mathscr{L},\Rightarrow ),\vdash)$ be a normal sequent structure and suppose that for every $\mathsf{Seq}(\mathscr{L},\Rightarrow )$-provable sequent $\Gamma\Rightarrow\Delta$, the following statements hold. 
\begin{enumerate}[label=(\arabic*)]
    \item If $\Gamma$ is not the empty list, then there exists a rule $R$ and a sequent $\Gamma^R\Rightarrow\Delta^R$ that is deducible from a set of sequents which includes $\Gamma\Rightarrow \Delta$ as one of the premise(s), by an application of the rule $R$. If $\Delta$ is not the empty list, then there exists a rule $L$ and a sequent $\Gamma^L\Rightarrow\Delta^L$ such that a similar conclusion holds. Furthermore, each of the premises are provable.
    \item{\label{nseqsys(3)} If there exists a proof of $\Gamma^R\Rightarrow \Delta^R$, then there exists a $T$-free $\mathsf{Seq}(\mathscr{L},\Rightarrow)$-proof of $\Gamma\Rightarrow \Delta$. Similarly, for $\Gamma^L\Rightarrow \Delta^L$.}
    \end{enumerate}
Then every provable sequent has a $T$-free proof.
\end{thm}}

\begin{proof}
    Consider any $\mathsf{Seq}(\mathscr{L},\Rightarrow)$-provable $\Gamma\Rightarrow \Delta$. Without loss of generality, let us assume that $\Gamma$ is nonempty. Then, by (1), it follows that there exists a rule $R$ and a sequent $\Gamma^R\Rightarrow \Delta^R$ such that it is deducible from a set of sequents which includes $\Gamma\Rightarrow \Delta$ is one of the premises and $\Gamma^R\Rightarrow \Delta^R$ the conclusion, of some instance of $R$. Furthermore, since each of the premises is provable, it follows that $\Gamma^R\Rightarrow \Delta^R$ is provable as well. Consequently, by (2), it follows that $\Gamma\Rightarrow \Delta$ has a $T$-free proof.    
\end{proof}

\begin{thm}
Let $(\mathsf{Seq}(\mathscr{L},\Rightarrow ),\vdash)$ be a normal sequent structure. Then, the following statements are equivalent.
\begin{enumerate}[label=(\arabic*)]
    \item Every provable sequent has a $T$-free proof.
    \item The conditions of \hyperlink{nseqsyscutelim2}{\thmref{nseqsyscutelim2}} are satisfied. 
    \end{enumerate}
\end{thm}

\begin{proof}
    We only show the proof of $(1)\implies (2)$ because the other implication is precisely \hyperlink{nseqsyscutelim2}{\thmref{nseqsyscutelim2}}.

    Note that, since, by hypothesis $\Gamma\Rightarrow \Delta$ is provable, so are $A, \Gamma\Rightarrow \Delta$ (by ${\color{red}{\mathsf{WL}}}$) and $\Gamma\Rightarrow \Delta, A$ (by ${\color{red}{\mathsf{WR}}}$). We take $R={\color{red}{\mathsf{WR}}}$ and $L={\color{red}{\mathsf{WL}}}$. In other words, both $A, \Gamma\Rightarrow \Delta$ and $\Gamma\Rightarrow \Delta, A$ are deducible from $ \Gamma\Rightarrow \Delta$. Since $\Gamma\Rightarrow \Delta$ is provable, it follows that both $A, \Gamma\Rightarrow \Delta$ and $\Gamma\Rightarrow \Delta, A$ are provable. This proves (1) of \hyperlink{nseqsyscutelim2}{\thmref{nseqsyscutelim2}}. The proof of (2) of \hyperlink{nseqsyscutelim2}{\thmref{nseqsyscutelim2}} is immediate since every provable sequent has a $T$-free proof. Since $\Gamma\Rightarrow\Delta$ was chosen arbitrarily, we have thus shown (2).
\end{proof}

\begin{rem}
The proof of \hyperlink{nseqsyscutelim1}{\thmref{nseqsyscutelim1}} is an abstraction of the argument given in \hyperlink{ce1}{\thmref{ce1}}, and the proof of \hyperlink{nseqsyscutelim2}{\thmref{nseqsyscutelim2}} is an abstraction of the argument given in \hyperlink{thm:con+cut_elim(IV)}{\thmref{thm:con+cut_elim(IV)}}. In the proof of \hyperlink{thm:con+cut_elim(IV)}{\thmref{thm:con+cut_elim(IV)}}, the rules $R$ and $L$ were both ${\color{red}{\to\sf{L}}}$. 
\end{rem}

\begin{rem}
One may take the rule $T$ to be \textsc{cut}, or some restricted version of \textsc{cut} as needed. In general, there is no restriction on how $T$ is chosen, apart from the fact that it should be disjoint from both \textsc{weakening} and \textsc{exchange}.
\end{rem}
\begin{rem}
If a sequent structure has a witness function, then there may be more than one rule which it does not respect. 
\end{rem}

\section{Abstract Sequent Structures}
In the previous section, we have seen how to eliminate a rule $T$, provided we are dealing with normal sequent structures. Note that in the whole discussion we have not specified what the rule $T$ might be; but the sequents were constructed based on some algebraic language. In this section, we remove this restriction.

Note that since we are not considering an algebraic language, the elements of the sequents do not arise from an algebra. Hence, the rules that can potentially be eliminated via the methods discussed in this section must be purely “structural” in nature. 

In this section, we propose a generalisation of normal sequent structures, the latter connecting to the theory of
logical structures. We prove a \textsc{rule}-elimination theorem for a particular type of these structures and show that the converse of this theorem holds as well.

We start by analysing the proof of \hyperlink{nseqsyscutelim1}{\thmref{nseqsyscutelim1}} in the previous section. We observe the following.
\begin{enumerate}[label=(\alph*)]
\item As noted earlier, even though we used $\mathbb{N}$ in the above section, the only property of $\mathbb{N}$ used is that it is a well-ordered set. More specifically, we have relied on the \textit{inductive} structure of $\mathbb{N}$.
\item Even though we have assumed the notion of proof/deduction as is standard in the literature on sequent calculus, the properties that we required are limited only to results like \hyperlink{ml1}{\thmref{ml1}} and \hyperlink{ml2}{\thmref{ml2}} (or some suitable generalisation of it).
\item The base case of \hyperlink{nseqsyscutelim1}{\thmref{nseqsyscutelim1}} shows that if the value of $v$ on a given sequent is the least, then it has a $T$-free proof (where $T$ is disjoint from both \textsc{weakening} and \textsc{exchange}). 
\end{enumerate}

\begin{defn}[\textsc{Abstract Sequent Structure}] Let $\mathscr{L}$ be a set. An \textit{abstract sequent structure} is a tuple $\bf{AbSS}:=\bigl(\mathsf{Seq}(\mathscr{L},\Rightarrow ),\vdash\bigr)$ where,
\begin{enumerate}[label=(\roman*)]
\item $\mathsf{Seq}(\mathscr{L},\Rightarrow )$ denotes the set of all $\mathscr{L}$-sequents with respect to $\Rightarrow$ for some set $\mathscr{L}$.
\item $\vdash\,\subseteq\,\mathcal{P}(\mathsf{Seq}(\mathscr{L},\Rightarrow ))\times \mathsf{Seq}(\mathscr{L},\Rightarrow )$
\end{enumerate}\end{defn}

\begin{defn}[\textsc{A}b\textsc{SS Sequent}]Let $\mathscr{L}$ be a set and let $\Rightarrow\subseteq \boldsymbol{O}(\mathscr{L})\times\boldsymbol{O}(\mathscr{L})$, where $\boldsymbol{O}(\mathscr{L})$ denotes the class of all ordered lists of $\mathscr{L}$.%
\footnote{More formally, $\boldsymbol{O}(\mathscr{L}):=\{\mathscr{L}^I:I~\text{is a set}\}$; hence it is a class.} Then $(\Gamma,\Delta)\in \Rightarrow$ will be said to be a \textit{$\mathscr{L}$-sequent with respect to $\Rightarrow$}. In this case, instead of writing $(\Gamma,\Delta)\in \Rightarrow$, we will write $\Gamma\Rightarrow \Delta$. The \textit{set of all $\mathscr{L}$-sequents with respect to $\Rightarrow$} will henceforth be denoted by $\mathsf{Seq}(\mathscr{L},\Rightarrow )$.%
\end{defn}

\hypertarget{def:abss_rule}{\begin{defn}[\textsc{A}b\textsc{SS Rule}]{\label{def:abss_rule}}Let $\mathsf{Seq}(\mathscr{L},\Rightarrow )$ be the set of all $\mathscr{L}$-sequents with respect to $\Rightarrow$. A set $R$ is said to be a \textit{$\mathsf{Seq}(\mathscr{L},\Rightarrow )$-rule} if, $$\emptyset \subsetneq R\subseteq \mathsf{Seq}(\mathscr{L},\Rightarrow )^I\times \mathsf{Seq}(\mathscr{L},\Rightarrow )$$for some nonempty index set $I$. It is said to be a \textit{$\mathsf{Seq}(\mathscr{L},\Rightarrow )$-axiomatic rule} if $\emptyset \subsetneq R\subseteq \{\emptyset\}\times \mathsf{Seq}(\mathscr{L},\Rightarrow )$; otherwise it is called a \textit{non-axiomatic rule}. An \textit{instance} of a rule $R$ is an element of $R$. Given $((\Gamma_j\Rightarrow \Delta_j)_{j\in J} ,\Gamma\Rightarrow \Delta)\in R$, $\Gamma_i\Rightarrow \Delta_i$ is to be a \textit{premise} of this instance of $R$ for each $i\in \{1,\ldots,n\}$ and $\Gamma\Rightarrow \Delta$ is said to be a \textit{conclusion} of this instance of $R$. In this case, we will also sometimes say $\Gamma\Rightarrow \Delta$ to be \textit{obtained by an application of $R$}.\end{defn}}

 
\hypertarget{defn:RIAbSS}{\begin{defn}[\textsc{Rule-Induced A}b\textsc{SS}]{\label{def:RIAbSS}}Let $\bf{AbSS}:=\bigl(\mathsf{Seq}(\mathscr{L},\Rightarrow ),\deduc\bigr)$ be an abstract sequent structure and $\mathcal{R}:=\{R_i\}_{i\in I}$ be a family of $\bf{AbSS}$-rules. Then, $\bf{AbSS}:=\bigl(\mathsf{Seq}(\mathscr{L},\Rightarrow ),\deduc\bigr)$ is said to be \textit{induced/generated by $\{R_i\}_{i\in I}$} if the following statements hold.
\begin{enumerate}[label=(\roman*)]
\item $\emptyset\deduc \Gamma\Rightarrow  \Delta$ if $(\emptyset,\Gamma\Rightarrow \Delta)\in R_i$ whenever $R_i$ is an $\bf{AbSS}$-axiomatic rule.
\item $\emptyset\deduc \Gamma\Rightarrow  \Delta$ implies that there exists $R_i$ such that $(\emptyset,\Gamma\Rightarrow \Delta)\in R_i$ where $R_i$ is a $\bf{AbSS}$-axiomatic rule or  $((\Gamma_j\Rightarrow \Delta_j)_{j\in J},\Gamma\Rightarrow \Delta)\in R_i$ for some family of sequents $(\Gamma_j\Rightarrow \Delta_j)_{j\in J}$ such that $\emptyset\deduc_{\mathcal{R}} \Gamma_j\Rightarrow\Delta_j$ for all $j\in J$.
\item For all non-axiomatic $\bf{AbSS}$-rule $R_i$ and $((\Gamma_j \Rightarrow \Delta_j)_{j\in J}, \Gamma\Rightarrow\Delta) \in R_i$, the following statement holds: if $\emptyset\deduc_{\mathcal{R}} \Gamma_j\Rightarrow\Delta_j$ for each $j\in J$ then $\emptyset\deduc_{\mathcal{R}}\Gamma\Rightarrow\Delta$.
\end{enumerate}In this case, we will write $\deduc\,=\,\deduc_{\mathcal{R}}$.\end{defn}}

\begin{rem}
In informal terms, the first property asserts that in an $\mathbf{AbSS}$, axiomatic inferences are provable; the second one asserts that if there is a proof of some sequent then either it is the conclusion of an axiomatic rule or is deduced by the application of some non-axiomatic rule, each of whose premise sequents are provable; the last property says that if each premise of an instance of a rule is provable then so is its conclusion. 
\end{rem}

\begin{rem}
It is important to note that \hyperlink{def:RIAbSS}{\defref{def:RIAbSS}} is motivated by the `root-based' definition of proofs instead of the `leaf-based' one that we have used till \S3 in the context of \textbf{PLK}, or for normal sequent structures (for details and differences between these two approaches, see \cite[\S1.6]{Indrzejczak2021} and \cite{Segerberg1982}). Nevertheless, there is an important difference: we are not giving an explicit definition of proof or deduction; instead, we are simply specifying certain properties that an abstract sequent structure must satisfy in order to be a rule-induced abstract sequent structure. 
\end{rem}

\begin{defn}[\textsc{Rule-Induced A}b\textsc{SS-Eliminable Rule}]Given an abstract sequent structure $\bf{AbSS}:=\bigl(\mathsf{Seq}(\mathscr{L},\Rightarrow ),\vdash\bigr)$ generated by a family $\mathcal{R}:=\{R_i\}_{i\in I}$ of $\bf{AbSS}$-rule(s), an $\bf{AbSS}$-rule $R_i$ is said to be \textit{$\bf{AbSS}$-eliminable} if $\emptyset\vdash\Gamma\Rightarrow  \Delta$ implies that $\emptyset\vdash_{\mathcal{R}\setminus\{R_i\}} \Gamma\Rightarrow  \Delta$, i.e., if $\emptyset\vdash_{\mathcal{R}}\Gamma\Rightarrow  \Delta$ implies that $\emptyset\vdash_{\mathcal{R}\setminus\{R_i\}} \Gamma\Rightarrow  \Delta$.\end{defn}

\begin{defn}Let $\bf{AbSS}:=\bigl(\mathsf{Seq}(\mathscr{L},\Rightarrow ),\vdash\bigr)$ be an abstract sequent structure generated by a family $\mathcal{R}:=\{R_i\}_{i\in I}$ of $\bf{AbSS}$-rule(s). Let $(P,\le)$ be a poset and $v:\mathsf{Seq}(\mathscr{L},\Rightarrow )\to P$. Then $v$ is said to \textit{respect} an instance $((\Gamma_j\Rightarrow \Delta_j)_{j\in J},\Gamma\Rightarrow \Delta)\in R_i$ of the $\bf{AbSS}$-rule $R_i$ if $v(\Gamma_j\Rightarrow \Delta_j)\le v(\Gamma\Rightarrow \Delta)$ for all $j\in J$. It is said to  \textit{strictly respect} $((\Gamma_j\Rightarrow \Delta_j)_{j\in J},\Gamma\Rightarrow \Delta)$ if $v(\Gamma_j\Rightarrow \Delta_j)< v(\Gamma\Rightarrow \Delta)$ for all $j\in J$.\end{defn}

\hypertarget{rem:RILS}{\begin{rem}{\label{rem:RILS}}
Even though we have formulated the notion of rule-induced abstract sequent structures using the language of sequents, it is not necessary. Nevertheless, we have chosen to do so mainly due to the continuity of the presentation.  
\end{rem}}

\begin{defn}A \textit{relation space} is a pair $(P,\mathfrak{R})$ where $P$ is a set and $\mathfrak{R} \subseteq P\times P$.
\end{defn}

\hypertarget{ruleelim2}{\begin{thm}[\textsc{Rule-Elimination Theorem for A}b\textsc{SSs}]{\label{ruleelim2}}Let $\bf{AbSS}:=\bigl(\mathsf{Seq}(\mathscr{L},\Rightarrow ),\deduc_{\mathcal{R}}\bigr)$ be an abstract sequent structure generated by a family $\mathcal{R}:=\{R_i\}_{i\in I}$ of $\bf{AbSS}$-rules and let $(P,\mathfrak{R})$ be a relation space. Let $v:\mathsf{Seq}(\mathscr{L},\Rightarrow )\to P$ and $R_k$ be an $\bf{AbSS}$-non-axiomatic rule. Suppose that the following conditions are satisfied.
\begin{enumerate}[label=(\arabic*)]
\item $R_i\cap R_k=\emptyset$ for all $R_i\in \mathcal{R}\setminus\{R_k\}$
\item For every $((\Gamma_t\Rightarrow \Delta_t)_{t\in T},\Gamma\Rightarrow \Delta)\in R_k$, $(v(\Gamma_t\Rightarrow \Delta_t),v(\Gamma\Rightarrow \Delta))\in \mathfrak{R}$ for all $t\in T$.
\item If $\emptyset\deduc_{\mathcal{R}}\Gamma\Rightarrow\Delta$ then either $(\emptyset, \Gamma \Rightarrow \Delta)\in R_i$ for some $\bf{AbSS}$-axiomatic rule or there exists some $\bf{AbSS}$-non-axiomatic rule $R_i$ and $((\Gamma_j \Rightarrow \Delta_j)_{j\in J}, \Gamma\Rightarrow\Delta) \in R_i$ such that $\emptyset\deduc_{\mathcal{R}} \Gamma_j\Rightarrow\Delta_j$ for all $j\in J$ and $(v(\Gamma_j \Rightarrow \Delta_j), v(\Gamma\Rightarrow\Delta))\notin\mathfrak{R}$ for some $j \in J$. (Note that any $R_i$ satisfying the second condition should be disjoint from $R_k$ because of the previous condition.)
\item If $((\Gamma_j\Rightarrow \Delta_j)_{j\in J},\Gamma\Rightarrow \Delta)\in R_i$ for some $\bf{AbSS}$-rule $R_i$ such that $R_i\cap R_k=\emptyset$ with $(v(\Gamma_j\Rightarrow \Delta_j), v(\Gamma\Rightarrow \Delta))\notin \mathfrak{R}$ for some $j\in J$, and  $\emptyset\deduc_{\mathcal{R}}\Gamma_j\Rightarrow \Delta_j$ for all $j\in J$, then $\emptyset\deduc_{\mathcal{R}\setminus\{R_k\}}\Gamma_j\Rightarrow \Delta_j$ for all $j\in J$.
\end{enumerate}Then $R_k$ is $\bf{AbSS}$-eliminable.\end{thm}}

\begin{proof}Let $\emptyset\deduc_{\mathcal{R}}\Gamma\Rightarrow  \Delta$.
If $(\emptyset,\Gamma\Rightarrow  \Delta)\in R_i$ for some $\bf{AbSS}$-axiomatic rule, it follows that $\emptyset\deduc_{\mathcal{R}\setminus\{R_k\}}\Gamma\Rightarrow \Delta$. Otherwise, by (3), there exists a non-axiomatic $R_i$ with $((\Gamma_j\Rightarrow \Delta_j)_{j\in J},\Gamma\Rightarrow \Delta)\in R_i$ such that $R_i\cap R_k=\emptyset$, $\emptyset\deduc_{\mathcal{R}}\Gamma_j\Rightarrow \Delta_j$ for every $j\in J$ and $(v(\Gamma_j \Rightarrow \Delta_j), v(\Gamma\Rightarrow\Delta))\notin\mathfrak{R}$ for some $j \in J$. Then, by (4), it follows that $\emptyset\deduc_{\mathcal{R}\setminus\{R_k\}}\Gamma_j\Rightarrow  \Delta_j$ for all $j\in J$. Consequently, since $R_{i}\notin \mathcal{R}\setminus\{R_k\}$, it follows that, $\emptyset\deduc_{\mathcal{R}\setminus\{R_k\}}\Gamma\Rightarrow  \Delta$. So, $R_k$ is $\bf{AbSS}$-eliminable.\end{proof}



\hypertarget{ruleelim2conv}{\begin{thm}{\label{ruleelim2conv}}Let $\bf{AbSS}:=\bigl(\mathsf{Seq}(\mathscr{L},\Rightarrow ),\deduc_{\mathcal{R}}\bigr)$ be an abstract sequent structure generated by a family $\mathcal{R}:=\{R_i\}_{i\in I}$ of $\bf{AbSS}$-rule(s). Let $R_k$ be an $\bf{AbSS}$-eliminable non-axiomatic rule such that it is disjoint from every rule in $\mathcal{R}\setminus\{R_k\}$. Then there exists a relation space $(P,\mathfrak{R})$ and a function $v:\mathsf{Seq}(\mathscr{L},\Rightarrow )\to P$ such that,
\begin{enumerate}[label=(\arabic*)]
\item For every $((\Gamma_t\Rightarrow \Delta_t)_{t\in T},\Gamma\Rightarrow \Delta)\in R_k$, $(v(\Gamma_j\Rightarrow \Delta_j), v(\Gamma\Rightarrow \Delta))\in\mathfrak{R}$ for all $t\in T$.
\item If $\emptyset\deduc_{\mathcal{R}}\Gamma\Rightarrow  \Delta$ then either $(\emptyset, \Gamma\Rightarrow  \Delta)\in R_i$ for some $\bf{AbSS}$-axiomatic rule or $((\Gamma_j\Rightarrow \Delta_j)_{j\in J},\Gamma\Rightarrow \Delta)\in R_i$ for some $\bf{AbSS}$-rule $R_i$ such that $\emptyset\deduc_{\mathcal{R}} \Gamma_j\Rightarrow\Delta_j$ for all $j\in J$ and $(v(\Gamma_j \Rightarrow \Delta_j), v(\Gamma\Rightarrow\Delta))\notin\mathfrak{R}$ for some $j \in J$.
\item If $((\Gamma_j\Rightarrow \Delta_j)_{j\in J},\Gamma\Rightarrow \Delta)\in R_i$ for some $\bf{AbSS}$-rule $R_i$ with $(v(\Gamma_j\Rightarrow \Delta_j), v(\Gamma\Rightarrow \Delta))\notin \mathfrak{R}$ and $\emptyset\deduc_{\mathcal{R}}\Gamma_j\Rightarrow \Delta_j$ for all $j\in J$ then $\emptyset\deduc_{\mathcal{R}\setminus\{R_k\}}\Gamma_j\Rightarrow \Delta_j$ for all $j\in J$.
\end{enumerate}\end{thm}}

\begin{proof}
We define a relation space $(\mathsf{Seq}(\mathscr{L},\Rightarrow ),\mathfrak{R})$ as follows,

\begin{align*}
    (\Gamma\Rightarrow \Delta, \Pi\Rightarrow \Lambda)\in \mathfrak{R}\iff&\Gamma\Rightarrow \Delta~ \text{is a premise of an instance of}~R_k~\text{and} \\&\Pi\Rightarrow \Lambda~\text{is the conclusion of the same instance}
\end{align*}The function $v$ is the identity function.

(1) is immediately satisfied by observing the relation we defined.

For (2), suppose $\emptyset\deduc_{\mathcal{R}} \Gamma\Rightarrow \Delta$. Since $R_k$ is eliminable, $\emptyset\deduc_{\mathcal{R}\setminus\{R_k\}}\Gamma\Rightarrow \Delta$. Consequently, since $\bf{AbSS}$ is an abstract sequent structure, it follows that either there exists $R_i$ such that $(\emptyset,\Gamma\Rightarrow \Delta)\in R_i$ where $R_i$ is an $\bf{AbSS}$-axiomatic rule (since $R_k$ is non-axiomatic and hence the sets of axiomatic rules in $\mathcal{R}$ and $\mathcal{R}\setminus\{R_k\}$ are same) or  $((\Gamma_j\Rightarrow \Delta_j)_{j\in J},\Gamma\Rightarrow \Delta)\in R_i\in \mathcal{R}\setminus\{R_k\}$ for some family of sequents $(\Gamma_j\Rightarrow \Delta_j)_{j\in J}$ such that $\emptyset\deduc_{\mathcal{R}\setminus\{R_k\}}\Gamma_j\Rightarrow\Delta_j$ for all $j\in J$. In the former case, the conclusion is immediate. Otherwise, since $R_i\in \mathcal{R}\setminus\{R_k\}$ and it is disjoint from $R_k$, the latter conclusion follows.

For (3), note that if $(\Gamma_j\Rightarrow \Delta_j, \Gamma\Rightarrow \Delta)\notin \mathfrak{R}$ for some $j\in J$ with $((\Gamma_j\Rightarrow \Delta_j)_{j\in J},\Gamma\Rightarrow \Delta)\in R_i$ for some $\bf{AbSS}$-rule $R_i$, it follows that $R_i\ne R_k$. Hence, by hypothesis, $R_i\cap R_k=\emptyset$. So $R_i\in \mathcal{R}\setminus\{R_k\}$. Now, if $((\Gamma_j\Rightarrow \Delta_j)_{j\in J},\Gamma\Rightarrow \Delta)\in R_i$ for some $\bf{AbSS}$-rule $R_i\in \mathcal{R}\setminus\{R_k\}$ and if $\emptyset\deduc_{\mathcal{R}}\Gamma_j\Rightarrow \Delta_j$ for all $j\in J$ then since $R_k$ is eliminable, we immediately have $\emptyset\deduc_{\mathcal{R}\setminus\{R_k\}}\Gamma_j\Rightarrow \Delta_j$ for all $j\in J$. This completes the proof.
\end{proof}

\begin{rem}
The condition of mutual disjointness of the family of rules in the above theorems is due to the intuition that instances of two distinct rules should be \textit{formally/schematically} distinguishable. 
\end{rem}

\section{Concluding Remarks}
In this paper, we have proved several \textsc{rule}-elimination theorems. These theorems are intended to provide insight into the basic structure of the proofs of the elimination of `structural rules' from a sequent system. Even though \textsc{cut}-elimination can be seen as a special case of these theorems, applying these theorems is difficult due mainly to the level of abstraction. Finding some more easily applicable criteria remains a project for the future. 

Even though our discussion centred only around the single-conclusion sequent structures, we think that a similar theory can be worked out for multiple conclusion sequent structures with suitable modifications. One of our long-term aims is to capture in a single theoretical framework most of the generalisations of sequent calculus and prove \textsc{rule}-elimination theorems for them. And, it seems to us that the problem can be stated even more generally independent of the notion of `sequents'. More specifically, as we have hinted in \hyperlink{rem:RILS}{\remref{rem:RILS}}, one can define something like `rule-induced logical structure' and investigate similar questions; thus connecting directly to the theory of logical structures. Satisfactory answers to such questions will not only have important consequences for proof-theoretic studies but also for the model-theoretic investigations as well.  

The idea of using an \textit{arbitrary} poset as the value set of a witness function can be useful, in particular, if one wants to carry out `induction-like' arguments. We would like to work with the ideas introduced in \cite{Ivanov2022} to see if it is possible to give an inductive version of the \textsc{rule}-elimination theorem for abstract sequent structures where the value set is an arbitrary poset.

\section*{Acknowledgements}

The primary motivation for this work came from Prof. Jean-Yves B{\'e}ziau's tutorial on Universal Logic at the 7th World Congress and School on Universal Logic in Crete, Greece. 

I would like to thank Iris van der Giessen and Prof. Dick de Jongh for kindly providing detailed feedback on an earlier version of this work. Discussions with  Prof. Rosalie Iemhoff were invaluable. The comments of Prof. Reuben Rowe and Prof. Lloyd Humberstone were extremely helpful and improved the presentation of this article. I would also like to thank Prof. Andrzej Indrzejczak and Prof. Jayanta Sen for clarifying several doubts regarding sequent calculus. 

I am grateful to IIIT-D for partially funding my trips to several conferences where this paper was presented. 

Finally, and most importantly, this work would have been impossible without the constant support, encouragement and constructive criticisms of my advisors: Prof. Sankha S. Basu and Prof. Mihir K. Chakraborty. I am indebted to both of them beyond measure.  

\bibliographystyle{plain}
\bibliography{cutelim}
\end{document}